\documentclass[12pt]{elsarticle}
\usepackage{lineno,hyperref}
\usepackage[margin=1.0in]{geometry}
\usepackage{amsmath, amssymb,epsfig,bm}
\usepackage{latexsym}
\usepackage{algorithmicx}
\usepackage[ruled,vlined,linesnumbered]{algorithm2e}
\usepackage{graphicx,subfigure}
\usepackage{color}
\usepackage{appendix}

\modulolinenumbers[5]

\numberwithin{equation}{section}
\input{tcilatex}
\begin{document}

\begin{frontmatter}

\title{Numerical Solution of the 3-D Travel Time Tomography Problem}
\author[mk]{Michael V. Klibanov}
\cortext[mycorrespondingauthor]{Corresponding author}
\ead{mklibanv@uncc.edu}
\author[jz]{Jingzhi Li\corref{mycorrespondingauthor}}
\ead{li.jz@sustech.edu.cn}
\author[wl]{Wenlong Zhang\corref{mycorrespondingauthor}}
\ead{zhangwl@sustech.edu.cn}
\address[mk]{Department of Mathematics and Statistics, University of North Carolina at
Charlotte, Charlotte, NC, 28223, USA}
\address[jz]{Department of Mathematics \& National Center for Applied Mathematics
Shenzhen \& SUSTech International Center for Mathematics, Southern
University of Science and Technology, Shenzhen 518055, P.~R.~China}
\address[wl]{Department of Mathematics, Southern University of Science and Technology (SUSTech), 1088 Xueyuan Boulevard, University Town of Shenzhen, Xili, Nanshan, Shenzhen, Guangdong Province, P.R.China. }

\begin{abstract}
The first numerical solution of the 3-D travel time tomography problem is
presented. The globally convergent convexification numerical method is
applied.
\end{abstract}
\begin{keyword}
 \\ travel time tomography, numerical solution in 3D,
coefficient inverse problem, Carleman estimate, convexification, global
convergence

\textit{\textbf{AMS subject classification:} 35R30, 65M32.} 
\end{keyword}
\end{frontmatter}



\section{Introduction}

\label{sec:1}

In this paper the first 3-D computational result is obtained for the Travel
Time Tomography Problem (TTTP). We apply a globally convergent, the
so-called convexification method. Two versions of the theory of this method
for the TTTP\ were developed in \cite{kin1,kin2}. We use the version of \cite%
{kin1}, which is also fully described in the book \cite[Chapter 11]{KL}.

All functions considered below are real valued ones. Below $\mathbf{x}%
=\left( x,y,z\right) $ denotes points in $\mathbb{R}^{3}.$ Let $\Omega
\subset \mathbb{R}^{3}$ be a bounded domain and let $S\subset \mathbb{R}%
^{3}, $ $S\cap \overline{\Omega }=\varnothing $ be a surface, on a part $L_{%
\text{src}}\subseteq S$ of which wave sources are located. Let $c\left( 
\mathbf{x}\right) $ be the speed of waves propagation, $c\left( \mathbf{x}%
\right) =1/n\left( \mathbf{x}\right) $, where $n\left( \mathbf{x}\right) $
is the refractive index. Let $m\left( \mathbf{x}\right) =n^{2}\left( \mathbf{%
x}\right) .$ The function $m\left( \mathbf{x}\right) $ generates the
Riemannian metric \cite[Chapter 3]{Rom}%
\begin{equation*}
d\tau =\sqrt{m\left( \mathbf{x}\right) }\sqrt{\left( dx\right) ^{2}+\left(
dy\right) ^{2}+\left( dz\right) ^{2}}.
\end{equation*}%
The time, which the wave needs to propagate from a source $\mathbf{x}_{0}\in
L_{\text{src}}$ to a point $\mathbf{x}$ is called the the \textquotedblleft
first arrival time" or \textquotedblleft travel time". In the case of a
heterogeneous medium with $m\left( \mathbf{x}\right) \neq const.,$ the first
arriving signal, which arrives at the arrival time to the point $\mathbf{x}$%
, propagates not along the straight line connecting points $\mathbf{x}$ and $%
\mathbf{x}_{0}$ but rather along the geodesic line $\Gamma \left( \mathbf{x},%
\mathbf{x}_{0}\right) $ generated by this metric and connecting points $%
\mathbf{x}$ and $\mathbf{x}_{0}$. The travel time $\tau \left( \mathbf{x},%
\mathbf{x}_{0}\right) $ is%
\begin{equation}
\tau \left( \mathbf{x},\mathbf{x}_{0}\right) =\dint\limits_{\Gamma \left( 
\mathbf{x},\mathbf{x}_{0}\right) }\sqrt{m\left( \mathbf{y}\left( s\right)
\right) }ds,  \label{1.1}
\end{equation}%
where $ds$ is the element of the euclidean length. The function $\tau \left( 
\mathbf{x},\mathbf{x}_{0}\right) $ satisfies the so-called eikonal equation 
\cite[Chapter 3]{Rom}%
\begin{equation}
\left( \nabla _{\mathbf{x}}\tau \right) ^{2}=m\left( \mathbf{x}\right) ,
\label{1.2}
\end{equation}%
\begin{equation}
\tau \left( \mathbf{x},\mathbf{x}_{0}\right) =O\left( \left\vert \mathbf{x}-%
\mathbf{x}_{0}\right\vert \right) \text{ as }\mathbf{x}\rightarrow \mathbf{x}%
_{0}.  \label{1.20}
\end{equation}%
The Travel Time Tomography Problem (TTTP) is one of Coefficient Inverse
Problems (CIPs) for PDEs. Another name for TTTP is Inverse Kinematic Problem 
\cite[Chapter 3]{Rom}.

\textbf{Travel Time Tomography Problem (TTTP)}. \emph{Suppose that the
function }$m\left( \mathbf{x}\right) $\emph{\ is given outside of the domain 
}$\Omega .$\emph{\ Given boundary measurements }$g\left( \mathbf{x},\mathbf{x%
}_{0}\right) $\emph{\ of the function }$\tau \left( \mathbf{x},\mathbf{x}%
_{0}\right) ,$%
\begin{equation}
\tau \left( \mathbf{x},\mathbf{x}_{0}\right) =g\left( \mathbf{x},\mathbf{x}%
_{0}\right) ,\text{ }\mathbf{x}\in \partial \Omega ,\mathbf{x}_{0}\in L_{%
\text{src}},  \label{1.3}
\end{equation}%
\emph{find the function }$m\left( \mathbf{x}\right) $\emph{\ for }$\mathbf{x}%
\in \Omega .$

Thus, the governing eikonal PDE (\ref{1.2}) with condition (\ref{1.20}) is
nonlinear, its right hand side is unknown and geodesic lines $\Gamma \left( 
\mathbf{x},\mathbf{x}_{0}\right) $ in (\ref{1.1}) are unknown as well. These
three factors cause quite substantial challenges in attempts to solve this
problem.

TTTP has well known applications in the problem of the recovery of the speed
of propagation of seismic waves inside the Earth \cite{Herg,Rom,Vol,W}. Some
other applications are in the phaseless inverse scattering problem \cite{KR}%
, in the problem of detection and identification of underwater objects, and
in the problem of standoff inspection of buildings if using transmitted time
resolved electromagnetic data \cite{SAR1}. The first solution of TTTP was
obtained by Herglotz in 1905 \cite{Herg} and then by Wiechert and Zoeppritz 
\cite{W} in 1907. In these pioneering works the above application to
geophysics was considered, and the underlying mathematical model was 1-D. \
A detailed description of the method of \cite{Herg,W} can be found in \cite[%
section 3 of Chapter 3]{Rom}.

The $n-$D case, $n=2,3,$ is mathematically far more challenging than the 1-D
case. We refer to \cite{Sch}, where a numerical method was developed and
tested in the 2-D case. In \cite{Z} another numerical method was developed
analytically in the $n-$D, $n=2,3$ case and tested numerically in the 2$-$D
case. In \cite{kinlin} a numerical method with a guaranteed convergence was
developed for the linearized TTTP in the $n-$D case, assuming that the
function $m\left( \mathbf{x}\right) =m_{0}\left( \mathbf{x}\right)
+m_{1}\left( \mathbf{x}\right) ,$ where $m_{0}\left( \mathbf{x}\right) $ is
known, $m_{1}\left( \mathbf{x}\right) $ is unknown and $\left\vert
m_{1}\left( \mathbf{x}\right) \right\vert <<m_{0}\left( \mathbf{x}\right) ,$
in which case geodesic lines were known and generated by the function $%
m_{0}\left( \mathbf{x}\right) ,$ and the problem became linearized, see \cite%
[Chapter 3]{Rom} for a discussion of the importance of the linearized case.
In \cite{kinlin}, numerical studies were carried out in the 2$-$D case.

The authors are unaware about numerical results for TTTP in the full
nonlinear 3-D case, which would be supplied by the rigorous global
convergence analysis. The goal of the current paper is to obtain such
numerical results. As mentioned above, we implement numerically here the
globally convergent convexification method for the 3-D TTTP. The theory of
this method was published in \cite{kin1} and \cite[Chapter 11]{KL}. Our data
are formally determined ones. This means that the number $k$ of free
variables in the data equals the number $n$ of free variables in the unknown
function $m\left( \mathbf{x}\right) ,k=n=3.$ Also, our data are incomplete.
Indeed, in the case of complete data, the source should run somewhat around
the domain $\Omega .$ Unlike this, in our case $L_{\text{src}}\in \mathbb{R}%
^{3}\diagdown \Omega $ is an interval of a straight line.

All CIPs are both nonlinear and ill-posed. Conventional numerical methods
for CIPs are based on minimizations of least squares cost functionals, see,
e.g. \cite{Chavent,Gonch1,Gonch2}. However, nonlinearity and ill-posedness
of CIPs cause non convexity of these functionals, which, in turn typically
leads to the phenomenon of multiple local minima and ravines, see, e.g. \cite%
{Scales} for a numerical example of this phenomenon. However, since any
gradient-like method of the optimization of that functional can stop at any
point of a local minimum, which might be located far from the solution, then
it is worth to address that phenomenon. This is done by the convexification
concept.

The convexification concept was originally proposed in \cite{KI,Klib97} to
avoid the above phenomenon. Results of \cite{KI,Klib97} are only analytical
ones. Active numerical studies of various of versions of this concept have
started from the work \cite{Bak}, which has removed some obstacles for
numerical implementations. We refer here to, e.g. \cite%
{Khoa,KEIT,KLZ,KL,SAR2,TR} as some samples of these publications, also, see
references cited therein as well as the recent book \cite{KL}. The
convexification significantly modifies the idea of the paper \cite{BukhKlib}%
, in which the tool of Carleman estimates was introduced in the field of
inverse problems for the first time, see, e.g. books \cite{BK,BY,Is,KL,LRS}
for Carleman estimates. The idea of \cite{BukhKlib} has generated a number
of publications of many authors, which have discussed questions of
uniqueness and stability results for CIPs, see, e.g. \cite%
{BK,BY,Is,Klib84,Klib92,Ksurvey,KL,Yam} and references cited therein. Thus,
convexification stands aside of those publications, since it is dedicated to
the numerical extension of the idea of \cite{BukhKlib}. The numerical issue
was not discussed in \cite{BukhKlib}.

Given a CIP for a PDE, the convexification constructs a weighted
Tikhonov-like least squares cost functional $J_{\lambda }$ with the Carleman
Weight Function (CWF) in it. This is the function which is used as the
weight in the Carleman estimate for that PDE operator. Here $\lambda \geq 1$
is the parameter of the CWF. This functional is considered on a certain
convex bounded set $S\left( d\right) \subset H,$ where $d$ is the diameter
of this set and $H$ is an appropriate Hilbert space. The central theorem for
each version of the convexification is the one, which claims that if $%
\lambda $ is sufficiently large, then $J_{\lambda }$ is strictly convex on $%
S\left( d\right) $ and has unique minimizer on $\overline{S\left( d\right) }%
. $ Next, a theorem is proven, which claims convergence to that unique
minimizer of either the gradient projection or the gradient descent method
of the minimization of\ $J_{\lambda }$ with an arbitrary starting point of $%
S\left( d\right) .$ Since smallness restrictions are not imposed on $d$,
then this is the \emph{global} convergence. More precisely, we call a
numerical method for a CIP \emph{globally convergent} if a theorem is
proven, which claims that this method provides at least one point in a
sufficiently small neighborhood of the true solution without any advanced
knowledge of this neighborhood.

\textbf{Remarks 1.1:}

\begin{enumerate}
\item \emph{In the current paper, so as in \cite{kin1}, \cite[Chapter 11]{KL}%
, the CWF is generated by a Volterra-like integral operator. This is unlike
the above cited works on the convexification, in which CWFs are those
generated by PDE operators.}

\item \emph{Even though the theory of the convexification requires
sufficiently large values of the parameter }$\lambda ,$\emph{\ our numerical
experiments in, e.g. \cite{Khoa,KEIT,KLZ,KL,SAR2,TR} and many other related
publications on the convexification as well as in the current paper
consistently demonstrate that the choice }$\lambda \in \left[ 1,4\right] $%
\emph{\ is sufficient. Philosophically, this situation is similar with the
situation in any asymptotic theory. Such a theory usually states that if a
parameter }$X$\emph{\ is sufficiently large/small, then a certain formula }$%
Y $\emph{\ is valid with a good accuracy. However, in any specific problem
at hands with its specific choice of parameters only numerical experiments
can determine which exactly values of }$X$\emph{\ provide a good accuracy
for }$Y $\emph{. }
\end{enumerate}

In our approach, we use a certain approximate mathematical model. Our model
includes two elements: the truncation of a Fourier-like series with respect
to a special orthonormal basis\ in $L_{2}\left( a,b\right) $ \cite{KJIIP}, 
\cite[Chapter 6]{KL} and the so-called \textquotedblleft partial finite
differences" approach. In this approach finite differences are assumed with
respect to two spatial variables whereas the third one is treated in the
conventional continuos manner. The step size $h$ of these finite differences
is not \textquotedblleft allowed" to tend to zero. We do not know how to
prove convergence of our method neither in the case when the number of terms
of that series $N\rightarrow \infty $ nor in the case $h\rightarrow 0.$
Therefore, the \emph{only way} to verify the validity of this approximate
mathematical model is via numerical studies. In this regard, we note that
truncations of Fourier-like series with respect to the same basis were done
for a variety of CIPs in \cite{Khoa,KEIT,kinlin,TR} and \cite[Chapters
7,10,11,12]{KL}. In each of these cases, the validity of the corresponding
approximate mathematical model was verified numerically. Partial finite
differences were used in the numerical works \cite{Khoa,kinlin,TR}, where
other versions of the convexification were presented. Similar cases of
truncated Fourier series without proofs of convergence at $N\rightarrow
\infty $ can be observed for inverse problems considered by some other
authors, see, e.g. \cite{GN,Kab}. And successful numerical verifications
also took place in these references.

Philosophically, a similar situation with an approximate mathematical model
arises in optics since the Huygens-Fresnel theory is not yet rigorously
derived from the Maxwell's equations. Nevertheless, the Huygens-Fresnel
theory, which represents an approximate mathematical model of the
diffraction in optics, works quite well in practice and is, therefore,
commonly acceptable in optics. The corresponding discussion can be found in
section 8.1 of the classical textbook of Born and Wolf \cite{BW}.

In section 2 we construct the above mentioned globally strictly convex cost
functional for our CIP. In section 3 we formulate Theorems 3.1-3.4 of the
convergence analysis. Theorems 3.1 and 3.2 are known from \cite{kin1}, \cite[%
Chapter 11]{KL}. Theorems 3.3 and 3.4 are new ones. Theorem 3.3 is a new
accuracy estimate of the minimizer of our Tikhonov-like weighted cost
functional with the CWF in it. Theorem 3.4 claims the global convergence of
the gradient descent method. On the other hand, a more complicated to
implement gradient projection method was used in the convergence analysis of 
\cite{kin1}, \cite[Chapter 11]{KL}.Therefore, we prove in section 4 only
theorems 3.3 and 3.4. In section 5 we describe our numerical studies.

\section{A Globally Strictly Convex Cost Functional}

\label{sec:2}

We now specify the domain $\Omega $ and the line of sources $L_{\text{src}}.$
Consider two numbers $B,\rho >0$. We define the domain $\Omega \subset 
\mathbb{R}^{3}$ as%
\begin{equation}
\Omega =\left\{ \mathbf{x}=\left( x,y,z\right) :x,y\in \left( 0,1\right)
,z\in \left( B,B+\rho \right) \right\} .  \label{2.1}
\end{equation}%
The boundary $\partial \Omega $ of the domain $\Omega $ consists of three
parts,%
\begin{equation}
\partial \Omega =D_{B}\cup D_{B+\rho }\cup \Gamma ,  \label{2.2}
\end{equation}%
\begin{equation}
D_{B}=\left\{ \mathbf{x}=\left( x,y,z\right) :x,y\in \left( 0,1\right)
,z=B\right\} ,  \label{2.3}
\end{equation}%
\begin{equation}
D_{B+\rho }=\left\{ \mathbf{x}=\left( x,y,z\right) :x,y\in \left( 0,1\right)
,z=B+\rho \right\} ,  \label{2.4}
\end{equation}%
\begin{equation}
\Gamma =\partial \Omega \diagdown \left( D_{B}\cup D_{B+\rho }\right) .
\label{2.5}
\end{equation}%
Let $a,b,d,z_{0}$ be three numbers, where $a<b,z_{0}<B.$ We assume that the
source runs along the line $L_{\text{src}}$,%
\begin{equation}
L_{\text{src}}=\left\{ \mathbf{x=}\left( x,y,z\right) :x=\alpha \in \left(
a,b\right) ,y=d,z=z_{0}\right\} \cap \overline{\Omega }=\varnothing .
\label{2.8}
\end{equation}%
Hence, the source $\mathbf{x}_{\alpha }=\left( \alpha ,d,z_{0}\right) \in L_{%
\text{src}},\alpha \in \left( a,b\right) $.\ Below $\tau \left( \mathbf{x}%
,\alpha \right) $ means $\tau \left( \mathbf{x},\mathbf{x}_{\alpha }\right) $
with $\mathbf{x}_{\alpha }\in L_{\text{src}}.$

We impose the following conditions on the function $m\left( \mathbf{x}%
\right) :$%
\begin{equation}
m\left( \mathbf{x}\right) \geq 1,\quad \mathbf{x}\in \mathbb{R}^{3},
\label{2.9}
\end{equation}%
\begin{equation}
m\left( \mathbf{x}\right) =1,\quad \mathbf{x}\in \left\{ z<B\right\} \cup
\left\{ \left( x,y\right) \notin \left( 0,1\right) \times \left( 0,1\right)
\right\} ,  \label{2.10}
\end{equation}%
\begin{equation}
m\in C^{2}\left( \mathbb{R}^{3}\right) ,  \label{2.11}
\end{equation}%
\begin{equation}
m_{z}\left( \mathbf{x}\right) \geq 0\text{, }\mathbf{x}\in \overline{\Omega }%
.  \label{2.12}
\end{equation}

\textbf{Remark 2.1}. \emph{The monotonicity condition (\ref{2.12}) can also
be found in \cite[section 2 of Chapter 3]{Rom} and \cite{Vol}. Also a
similar condition was imposed in the originating works of 1905, 1907 \cite%
{Herg,W}, see \cite[section 3 of Chapter 3]{Rom}. }

Lemma 2.1 follows immediately from Lemma 4.1 of \cite{kin1} and Lemma 11.4.1
of \cite{KL}.

\textbf{Lemma 2.1.} \emph{There exists a constant }$c_{0}=c_{0}\left(
a,b,d,z_{0},B\right) >0$\emph{\ depending only on listed parameters such
that }%
\begin{equation}
\tau _{z}\left( \mathbf{x},\alpha \right) \geq c_{0},\forall \mathbf{x}\in 
\overline{\Omega },\forall \alpha \in \left[ a,b\right] .  \label{2.120}
\end{equation}

We assume below that the function\emph{\ }$\tau \left( \mathbf{x},\alpha
\right) \in C^{1}\left( \mathbb{R}^{3}\times \left[ a,b\right] \right) .$ We
impose the assumption of the regularity of geodesic lines:

\textbf{Regularity of Geodesic Lines}. \emph{Let} \emph{the pair of points }$%
\left( \mathbf{x},\mathbf{x}_{\alpha }\right) $\emph{\ }$\in \overline{%
\Omega }\times L_{\text{src}}.$ \emph{Then there exists unique geodesic line 
}$\Gamma \left( \mathbf{x},\alpha \right) $\emph{\ connecting these two
points and }$\Gamma \left( \mathbf{x},\alpha \right) \cap D_{B}\neq
\varnothing $\emph{.} \emph{Also, if a geodesic line, which starts at the
point }$\mathbf{x}_{\alpha }\in L_{\text{src}},$ \emph{intersects }$D_{B},$ 
\emph{then it intersects it at a single point. Next, that geodesic line
intersects }$\partial \Omega \diagdown D_{B}$\emph{\ at another single
point, see (\ref{2.2})-(\ref{2.5}). In addition, after intersecting }$%
\partial \Omega \diagdown D_{B},$\emph{\ this line goes away from }$\Omega .$%
\emph{\ In other words, this line is not reflected back from any point of
its intersection with }$\partial \Omega $.

It follows from (\ref{2.1}), (\ref{2.8}) and (\ref{2.10}) that 
\begin{equation*}
\tau \left( x,y,z,\alpha \right) =\sqrt{\left( x-\alpha \right) ^{2}+\left(
y-d\right) ^{2}+\left( z-z_{0}\right) ^{2}}\text{ for }z<B.
\end{equation*}%
Denote $u_{0}\left( \mathbf{x},\alpha \right) =\tau _{z}^{2}\left(
x,y,B,\alpha \right) $. Then 
\begin{equation}
u_{0}\left( \mathbf{x},\alpha \right) =\frac{\left( B-z_{0}\right) ^{2}}{%
\left( x-\alpha \right) ^{2}+\left( y-d\right) ^{2}+\left( B-z_{0}\right)
^{2}}.  \label{2.13}
\end{equation}

\subsection{A boundary value problem for a nonlinear integral differential
equation}

\label{sec:2.1}

Denote 
\begin{equation}
u\left( \mathbf{x},\alpha \right) =\tau _{z}^{2}\left( \mathbf{x},\alpha
\right) ,\mathbf{x}\in \Omega ,\alpha \in \left( a,b\right) .  \label{2.14}
\end{equation}%
By Lemma 2.1 
\begin{equation}
\tau _{z}\left( \mathbf{x},\alpha \right) =\sqrt{u\left( \mathbf{x},\alpha
\right) }\text{, }\mathbf{x}\in \Omega ,\alpha \in \left( a,b\right) .
\label{2.15}
\end{equation}%
Hence, (\ref{1.3}) and (\ref{2.15}) imply 
\begin{equation*}
\tau \left( \mathbf{x},\alpha \right) =-\dint\limits_{z}^{B+\rho }\sqrt{%
u\left( x,y,t,\alpha \right) }dt+g\left( x,y,B+\rho ,\alpha \right) ,\mathbf{%
x}\in \Omega ,\alpha \in \left( a,b\right) ,
\end{equation*}%
\begin{equation}
\tau _{x}\left( \mathbf{x},\alpha \right) =-\dint\limits_{z}^{B+\rho }\left( 
\frac{u_{x}}{2\sqrt{u}}\right) \left( x,y,t,\alpha \right) dt+g_{x}\left(
x,y,B+\rho ,\alpha \right) ,\mathbf{x}\in \Omega ,\alpha \in \left(
a,b\right) ,  \label{2.16}
\end{equation}%
\begin{equation}
\tau _{y}\left( \mathbf{x},\alpha \right) =-\dint\limits_{z}^{B+\rho }\left( 
\frac{u_{y}}{2\sqrt{u}}\right) \left( x,y,t,\alpha \right) dt+g_{y}\left(
x,y,B+\rho ,\alpha \right) ,\mathbf{x}\in \Omega ,\alpha \in \left(
a,b\right) .  \label{2.17}
\end{equation}%
Substituting (\ref{2.15})-(\ref{2.17}) in the eikonal equation (\ref{1.2}),
we obtain the following equation for $\mathbf{x}\in \Omega ,\alpha \in
\left( a,b\right) $: 
\begin{equation*}
u\left( \mathbf{x},\alpha \right) +\left[ -\int \displaylimits_{z}^{B+\rho
}\left( \frac{u_{x}}{2\sqrt{u}}\right) \left( x,y,t,\alpha \right)
dt+g_{x}\left( x,y,B+\rho ,\alpha \right) \right] ^{2}
\end{equation*}%
\begin{equation}
+\left[ -\int \displaylimits_{z}^{B+\rho }\left( \frac{u_{y}}{2\sqrt{u}}%
\right) \left( x,y,t,\alpha \right) dt+g_{y}\left( x,y,B+\rho ,\alpha
\right) \right] ^{2}=m\left( \mathbf{x}\right) .  \label{2.160}
\end{equation}%
Differentiating (\ref{2.160}) with respect to $\alpha $ and using $\partial
_{\alpha }m\left( \mathbf{x}\right) \equiv 0,$ we obtain 
\begin{equation}
u_{\alpha }\left( x,y,z,\alpha \right) +\frac{\partial }{\partial \alpha }%
\left[ -\dint\limits_{z}^{B+\rho }\left( \frac{u_{x}}{2\sqrt{u}}\right)
\left( x,y,t,\alpha \right) dt+g_{x}\left( x,y,B+\rho ,\alpha \right) \right]
^{2}  \label{2.18}
\end{equation}%
\begin{equation*}
+\frac{\partial }{\partial \alpha }\left[ -\dint\limits_{z}^{B+\rho }\left( 
\frac{u_{y}}{2\sqrt{u}}\right) \left( x,y,t,\alpha \right) dt+g_{y}\left(
x,y,B+\rho ,\alpha \right) \right] ^{2},\mathbf{x}\in \Omega ,\alpha \in
\left( a,b\right) .
\end{equation*}

Thus, we came up with the following boundary value problem (BVP).

\textbf{Boundary Value Problem 1 (BVP1).} \emph{Find the function }$u\left( 
\mathbf{x},\alpha \right) \in C^{1}\left( \overline{\Omega }\times \left[ a,b%
\right] \right) $\emph{\ satisfying both integral differential equation (\ref%
{2.18}) for }$\mathbf{x}\in \Omega ,\alpha \in \left( a,b\right) $\emph{\
and the boundary condition}%
\begin{equation}
u\mid _{\Gamma \cup D_{B}}=\widetilde{g}\left( \mathbf{x},\alpha \right) ,
\label{2.19}
\end{equation}%
\begin{equation}
\widetilde{g}\left( \mathbf{x},\alpha \right) =\left\{ 
\begin{array}{c}
g_{z}^{2}\left( \mathbf{x},\alpha \right) ,\text{ }\mathbf{x}\in \Gamma
,\alpha \in \left( a,b\right) , \\ 
u_{0}\left( \mathbf{x},\alpha \right) ,\text{ }\mathbf{x}\in D_{B},\alpha
\in \left( a,b\right)%
\end{array}%
\right.  \label{2.190}
\end{equation}%
\emph{where the function }$g\left( \mathbf{x},\alpha \right) $\emph{\ is
given in (\ref{1.3}) and the function }$u_{0}\left( \mathbf{x},\alpha
\right) $\emph{\ is given in (\ref{2.13}).}

Suppose that we have solved this problem. Then we substitute its solution in
the left hand side of equation (\ref{2.160}) and find the target function $%
m\left( \mathbf{x}\right) $. Clearly BVP1\ is a very complicated one.
Therefore, we construct below an approximate mathematical model for its
solution, and confirm the validity of this model computationally, see
section 1 for a relevant discussion.

\subsection{A special orthonormal basis in $L_{2}\left( a,b\right) $}

\label{sec:2.2}

This basis was first introduced in \cite{KJIIP} and then it was applied to a
number of other versions of the convexification method, see, e.g. \cite%
{Khoa,KEIT,kin1,kin2,kinlin,KL,TR}. Consider the set of functions $\left\{
q_{n}\left( \alpha \right) \right\} _{n=0}^{\infty }=\{\alpha ^{n}e^{\alpha
}\}_{n=0}^{\infty }\subset L_{2}(a,b).$ This is a set of linearly
independent functions, which is complete in the space $L_{2}(a,b)$. Applying
the Gram-Schmidt orthonormalization procedure to this set, we obtain the
orthonormal basis $\{\Phi _{n}\left( \alpha \right) \}_{n=0}^{\infty }$ in $%
L_{2}(a,b)$. Obviously, the function $\Phi _{n}(\alpha )$ has the form $\Phi
_{n}(\alpha )=Q_{n}(\alpha )e^{\alpha }$, $\forall n\geq 0,$ where $%
Q_{n}(\alpha )$ is a polynomial of the degree $n$. Let $a_{mn}=\left\{ \Phi
_{m},\Phi _{n}^{\prime }\right\} ,$ where $\left\{ ,\right\} $ is the scalar
product in $L_{2}(a,b)$. We have \cite{KJIIP}, \cite[Theorem 6.2.1]{KL}: 
\begin{equation}
a_{mn}=\left\{ 
\begin{array}{c}
1\text{ if }m=n, \\ 
0\text{ if }m>n.%
\end{array}%
\right.  \label{2.200}
\end{equation}%
Let $N\geq 1$ be an integer. Consider the $N\times N$ matrix $M_{N}=\left(
a_{mn}\right) _{m,n=0}^{N}.$ By (\ref{2.200}) $\det M_{N}=1.$ Thus, the
matrix $M_{N}$ is invertible. Note that neither classical orthonormal
polynomials nor the basis of trigonometric functions do not provide a
corresponding invertible matrix $M_{N}$ since the first function in such
cases is an identical constant, meaning that the first column of that analog
of $M_{N}$ is formed only by zeros.

We assume that the function $u$ can be represented via the truncated
Fourier-like series, 
\begin{equation}
u\left( \mathbf{x},\alpha \right) =\dsum\limits_{n=0}^{N-1}u_{n}\left( 
\mathbf{x}\right) \Phi _{n}\left( \alpha \right) ,\text{ }\mathbf{x}\in
\Omega ,\alpha \in \left( a,b\right) ,  \label{2.21}
\end{equation}%
\begin{equation}
W\left( \mathbf{x}\right) =\left( u_{0},...,u_{N-1}\right) ^{T}\left( 
\mathbf{x}\right) ,\text{ }\mathbf{x}\in \Omega .  \label{2.22}
\end{equation}%
Thus, the vector function $W\left( \mathbf{x}\right) $ is unknown. We also
assume that 
\begin{equation}
g\left( x,y,B+\rho ,\alpha \right) =\dsum\limits_{n=0}^{N-1}g_{n}\left(
x,y,B+\rho \right) \Phi _{n}\left( \alpha \right) ,\left( x,y\right) \in
\left( 0,1\right) ^{2},\alpha \in \left( a,b\right) ,  \label{2.23}
\end{equation}%
\begin{equation}
\widetilde{g}\left( \mathbf{x},\alpha \right) =\dsum\limits_{n=0}^{N-1}%
\widetilde{g}_{n}\left( \mathbf{x}\right) \Phi _{n}\left( \alpha \right) ,%
\text{ }\mathbf{x}\in \Gamma \cup D_{B},\text{ }\alpha \in \left( a,b\right)
,  \label{2.24}
\end{equation}%
\begin{equation}
G\left( x,y,B+\rho \right) =\left( g_{0},...,g_{N-1}\right) ^{T}\left(
x,y,B+\rho \right) ,\text{ }\widetilde{G}\left( \mathbf{x}\right) =\left( 
\widetilde{g}_{0},...\widetilde{g}_{N-1}\right) ^{T}\left( \mathbf{x}\right)
,  \label{2.25}
\end{equation}%
where functions $g$ and $\widetilde{g}$ are given in (\ref{1.3}) and (\ref%
{2.190}) respectively.

Substitute (\ref{2.21})-(\ref{2.25}) in (\ref{2.18}), assuming that (\ref%
{2.18}) holds for the truncated series (\ref{2.21}). Next, multiply
sequentially the obtained equality by the functions $\Phi _{n}\left( \alpha
\right) ,n=0,...,N-1$ and integrate with respect to $\alpha \in \left(
a,b\right) .$ We obtain the boundary value problem for the nonlinear system
of integral differential equations with respect to the vector function $%
W\left( \mathbf{x}\right) ,$%
\begin{equation}
M_{N}W+P\left( W_{x},W_{y},G_{x},G_{y},\mathbf{x}\right) =0,\mathbf{x}\in
\Omega ,  \label{2.26}
\end{equation}%
\begin{equation}
W\mid _{\Gamma \cup D_{B}}=\widetilde{G}\left( \mathbf{x}\right) .
\label{2.27}
\end{equation}%
where the $n-$th component of the $N-$D vector function $P=\left(
P_{0},...,P_{N-1}\right) ^{T}$ has the form 
\begin{equation*}
P_{n}\left( W_{x},W_{y},G_{x},G_{y},\mathbf{x}\right) =
\end{equation*}%
\begin{equation}
=\dint\limits_{a}^{b}\Phi _{n}\left( \alpha \right) \left[
-\dint\limits_{z}^{B+\rho }\left( \frac{u_{x}}{2\sqrt{u}}\right) \left(
x,y,t,\alpha \right) dt+g_{x}\left( x,y,B+\rho ,\alpha \right) \right]
^{2}d\alpha +  \label{2.28}
\end{equation}%
\begin{equation*}
+\dint\limits_{a}^{b}\Phi _{n}\left( \alpha \right) \left[
-\dint\limits_{z}^{B+\rho }\left( \frac{u_{y}}{2\sqrt{u}}\right) \left(
x,y,t,\alpha \right) dt+g_{y}\left( x,y,B+\rho ,\alpha \right) \right]
^{2}d\alpha ,\text{ }\mathbf{x}\in \Omega ,
\end{equation*}%
where $n=0,...,N-1$ and the function $u\left( \mathbf{x},\alpha \right) $ is
given in (\ref{2.21}).

Let $V$ be a Banach space. Consider the direct product $V\times V\times
...\times V,$ $N$ times. Then for the norm $\left\Vert \cdot \right\Vert
_{S} $ in $V$ we denote $\left\Vert \cdot \right\Vert _{S_{N}}$ an obvious
generalization of that norm for the case of the this direct product.

Let $R>0$ be an arbitrary number. Keeping in mind Lemma 2.1, (\ref{2.14})
and (\ref{2.27}), define the set $K\left( R\right) $ of $N-$D vector
functions $W\left( \mathbf{x}\right) $ as%
\begin{equation}
K\left( R\right) =\left\{ W\in C_{N}^{1}\left( \overline{\Omega }\right)
:u\left( \mathbf{x},\alpha \right) \geq c_{0}^{2}\text{ in }\overline{\Omega 
},W\mid _{\Gamma \cup D_{B}}=\widetilde{G}\left( \mathbf{x}\right)
,\left\Vert W\right\Vert _{C_{N}^{1}\left( \overline{\Omega }\right)
}<R\right\} ,  \label{2.29}
\end{equation}%
where $u\left( \mathbf{x},\alpha \right) $ is linked with $W$ via (\ref{2.21}%
), (\ref{2.22}), and $c_{0}$ is the number from (\ref{2.120}).

Thus, we have replaced the original BVP1 with an approximate BVP2:

\textbf{Boundary Value Problem 2 (BVP2). }\emph{Find the vector function }$%
W\in \overline{K\left( R\right) }$\emph{\ satisfying equation (\ref{2.26})
and condition (\ref{2.28}).}

\subsection{Partial finite differences}

\label{sec:2.3}

We now rewrite BVP2 via partial finite differences. Let $k>1$ be an integer.
Consider two partitions of the interval $\left( 0,1\right) $, see (\ref{2.1}%
), with the grid step size $h>0$:%
\begin{equation*}
0=x_{0}<x_{1}<...<x_{k}=1,x_{i+1}-x_{i}=h,i=0,...,k-1,
\end{equation*}%
\begin{equation*}
0=y_{0}<y_{1}<...<y_{k}=1,y_{i+1}-y_{i}=h,i=0,...,k-1,
\end{equation*}%
\begin{equation}
h\in \left[ h_{0},1\right) ,h_{0}=const.>0,  \label{2.32}
\end{equation}%
where $h_{0}$ is a fixed number. Let $\Omega _{1}=\left\{ \left( x,y\right)
:\left( x,y\right) \in \left( 0,1\right) \times \left( 0,1\right) \right\} .$
Define the discrete subset $\Omega ^{h}\hspace{0.5em}$of the domain $\Omega $
as:%
\begin{equation}
\Omega _{1}^{h}=\left\{ \left( x_{i},y_{j}\right) \right\} _{\left(
i,j\right) =\left( 0,0\right) }^{\left( i,j\right) =\left( k,k\right) },%
\text{ }  \label{2.33}
\end{equation}%
\begin{equation}
\Omega ^{h}=\left\{ \left( x_{i},y_{j},z\right) :\left\{ x_{i},y_{j}\right\}
_{\left( i,j\right) =\left( 0,0\right) }^{\left( i,j\right) =\left(
k,k\right) }\in \Omega _{1}^{h},z\in \left( B,B+\rho \right) \right\} .
\label{2.34}
\end{equation}%
We denote $\mathbf{x}^{h}=\left\{ \left( x_{i},y_{j},z\right) :\left(
x_{i},y_{j}\right) \in \Omega _{1}^{h},z\in \left( B,B+\rho \right) \right\}
.$ Similarly with (\ref{2.33}), (\ref{2.34}), we denote semidiscrete analogs
of parts of the boundary $\partial \Omega $ as $D_{B}^{h},D_{B+\rho
}^{h},\Gamma ^{h},$ and $\partial \Omega ^{h}=D_{B}^{h}\cup D_{B+\rho
}^{h}\cup \Gamma ^{h},$ see (\ref{2.2})-(\ref{2.5}). Let the vector function 
$Q(\mathbf{x})\in C_{N}^{1}(\overline{\Omega })$. Denote%
\begin{equation*}
Q_{ij}^{h}\left( z\right) =Q(x_{i},y_{j},z),\hspace{0.5em}i,j=0,...,k;\text{ 
}Q^{h}(\mathbf{x}^{h})=\left\{ Q^{h}\left( x_{i},y_{j},z\right) \right\} ,%
\text{ }z\in \left( a,b\right) .
\end{equation*}%
Thus, $Q^{h}(\mathbf{x}^{h})$ is an $N-$D vector function of discrete
variables $\left( x_{i},y_{j}\right) \in \Omega _{1}^{h}$ and continuous
variable $z\in \left( B,B+\rho \right) .$ Observe that the boundary terms of
this vector function at $\Gamma ^{h}$ are $\left\{ Q_{0,j}^{h}\left(
z\right) \right\} \cup \left\{ Q_{k,j}^{h}\left( z\right) \right\} \cup
\left\{ Q_{i,0}^{h}\left( z\right) \right\} \cup \left\{ Q_{i,k}^{h}\left(
z\right) \right\} ,i,j=0,...,k.$

We define finite difference derivatives of $Q^{h}(\mathbf{x}^{h})$ with
respect to $x,y$ only at interior points of the domain $\Omega ^{h}$ as: 
\begin{equation}
\partial _{x}Q^{h}(\mathbf{x}^{h})=\left\{ \frac{Q^{h}\left(
x_{i+1},y_{j},z\right) -Q^{h}\left( x_{i-1}y_{j},z\right) }{2h}\right\}
_{\left( i,j\right) =\left( 1,1\right) }^{\left( i,j\right) =\left(
k-1,k-1\right) },  \label{2.37}
\end{equation}%
\begin{equation}
\partial _{y}Q^{h}(\mathbf{x}^{h})=\left\{ \frac{Q^{h}\left(
x_{i},y_{j+1},z\right) -Q^{h}\left( x_{i},y_{j-1},z\right) }{2h}\right\}
_{\left( i,j\right) =\left( 1,1\right) }^{\left( i,j\right) =\left(
k-1,k-1\right) }.  \label{2.38}
\end{equation}%
We now define the semidiscrete analogs of spaces $L_{2,N}\left( \Omega
\right) $ and $H_{N}^{2}\left( \Omega \right) ,$ 
\begin{equation}
L_{2,N}^{h}\left( \Omega ^{h}\right) =\left\{ Q^{h}(\mathbf{x}%
^{h}):\left\Vert Q^{h}(\mathbf{x}^{h})\right\Vert _{L_{2,N}^{h}\left( \Omega
^{h}\right) }^{2}=\dsum\limits_{\left( i,j\right) =\left( 0,0\right)
}^{\left( k,k\right) }\left\{ \dint\limits_{a}^{b}\left( Q_{ij}^{h}\left(
z\right) \right) ^{2}dz\right\} <\infty \right\} ,  \label{2.40}
\end{equation}%
\begin{equation}
H_{N}^{1,h}\left( \Omega ^{h}\right) =\left\{ 
\begin{array}{c}
Q^{h}(\mathbf{x}^{h}):\left\Vert Q^{h}\right\Vert _{H_{N}^{1,h}\left( \Omega
^{h}\right) }^{2}= \\ 
\left\Vert Q_{x}^{h}\right\Vert _{L_{2N}^{h}\left( \Omega ^{h}\right)
}^{2}+\left\Vert Q_{y}^{h}\right\Vert _{L_{2N}^{h}\left( \Omega ^{h}\right)
}^{2}+\left\Vert Q^{h}\right\Vert _{L_{2N}^{h}\left( \Omega ^{h}\right)
}^{2}<\infty%
\end{array}%
\right\} .  \label{2.400}
\end{equation}%
\begin{equation}
H_{0,N}^{1,h}\left( \Omega ^{h}\right) =\left\{ Q^{h}(\mathbf{x}^{h})\in
H_{N}^{1,h}\left( \Omega ^{h}\right) :Q^{h}(\mathbf{x}^{h})\mid _{\Gamma
^{h}\cup D_{B}^{h}}=0\right\} .  \label{2.401}
\end{equation}%
It follows from (\ref{2.32}), (\ref{2.40}) and (\ref{2.400}) that norms in
the spaces (\ref{2.40}) and (\ref{2.400}) are equivalent. Thus, the space $%
H_{N}^{1,h}\left( \Omega ^{h}\right) $ is introduced only to stress the
presence of finite difference derivatives $Q_{x}^{h},Q_{y}^{h}.$

The following formulas are semidiscrete analogs of formulas (\ref{2.21}), (%
\ref{2.22}), (\ref{2.26})-(\ref{2.28}): 
\begin{equation}
u^{h}(\mathbf{x}^{h}\mathbf{,}\alpha )=\dsum\limits_{n=0}^{N-1}u_{n}^{h}(%
\mathbf{x}^{h})\Phi _{n}(\alpha ),  \label{2.42}
\end{equation}%
\begin{equation}
W^{h}(\mathbf{x}^{h})=\left( u_{0}^{h},...,u_{N-1}^{h}\right) ^{T}(\mathbf{x}%
^{h}),  \label{2.43}
\end{equation}%
\begin{equation}
M_{N}W^{h}=P\left( W_{x}^{h},W_{y}^{h},g_{x}^{h},g_{y}^{h},\mathbf{x}%
^{h}\right) ,\mathbf{x}^{h}\in \Omega ^{h},  \label{2.44}
\end{equation}%
\begin{equation}
W^{h}\mid _{\Gamma ^{h}\cup D_{B}^{h}}=\widetilde{G}^{h}\left( \mathbf{x}%
\right) ,  \label{2.45}
\end{equation}%
\begin{equation*}
P=\left( P_{0},...,P_{N-1}\right) ^{T},\text{ }P_{n}\left(
W_{x}^{h},W_{y}^{h},g_{x}^{h},g_{y}^{h},x_{i},y_{j},z\right) =
\end{equation*}%
\begin{equation}
=\dint\limits_{a}^{b}\Phi _{n}\left( \alpha \right) \left[
-\dint\limits_{z}^{B+\rho }\left( \frac{u_{x}^{h}}{2\sqrt{u^{h}}}\right)
\left( x_{i},y_{j},t,\alpha \right) dt+g_{x}^{h}\left( x_{i},y_{j},B+\rho
,\alpha \right) \right] ^{2}d\alpha +  \label{2.46}
\end{equation}%
\begin{equation*}
+\dint\limits_{a}^{b}\Phi _{n}\left( \alpha \right) \left[
-\dint\limits_{z}^{B+\rho }\left( \frac{u_{y}^{h}}{2\sqrt{u^{h}}}\right)
\left( x_{i},y_{j},t,\alpha \right) dt+g_{y}^{h}\left( x_{i},y_{j},B+\rho
,\alpha \right) \right] ^{2}d\alpha ,x^{h}\in \Omega ^{h}.
\end{equation*}%
We introduce the following semidiscrete analog of the set $K\left( R\right) $
in (\ref{2.29}) 
\begin{equation}
K^{h}\left( R\right) =\left\{ 
\begin{array}{c}
W^{h}\in H_{N}^{1,h}\left( \Omega ^{h}\right) :u^{h}\left( \mathbf{x},\alpha
\right) \geq c_{0}^{2}\text{ in }\overline{\Omega }^{h}, \\ 
W^{h}\mid _{\Gamma ^{h}\cup D_{B}^{h}}=\widetilde{G}^{h}\left( \mathbf{x}%
\right) ,\left\Vert W^{h}\right\Vert _{H_{N}^{1,h}\left( \Omega ^{h}\right)
}<R%
\end{array}%
\right\} .  \label{2.47}
\end{equation}

Thus, rather than solving BVP1 and BVP2, we solve below BVP3, which is the
semidiscrete analog of BVP2.

\textbf{Boundary Value Problem 3 (BVP3). }\emph{Find the vector function }$%
W^{h}\in K^{h}\left( R\right) $\emph{\ satisfying equations (\ref{2.42})-(%
\ref{2.46}).}

\subsection{Globally strictly convex cost functional for BVP3}

\label{sec:2.4}

First, we introduce the CWF. It is used below with the goal to arrange a
sort of the domination of the term $M_{N}W^{h}$ over the rest in (\ref{2.44}%
). Note that traditionally CWFs are used in the convexification only in
differential rather than in integral operators, see, e.g. the above cited
works on the convexification. Let $\lambda >0$ be a parameter. Our CWF is: 
\begin{equation}
\varphi _{\lambda }\left( z\right) =e^{2\lambda z}.  \label{2.48}
\end{equation}%
The proof of Lemma 2.1 is the same as the proof of Lemma 8.1 of \cite{kin1}
as well as of Lemma 11.8.1 of \cite{KL}.

\textbf{Lemma 2.1}. \emph{Let the function }$\varphi _{\lambda }\left(
z\right) $\emph{\ be the one defined in (\ref{2.48}). The following estimate
holds for all }$\lambda >0$\emph{\ and for all functions }$q\in L_{1}\left(
B,B+\rho \right) :$ 
\begin{equation*}
\dint\limits_{B}^{B+\rho }\left( \dint\limits_{z}^{B+\rho }\left\vert
q\left( y\right) \right\vert dy\right) \varphi _{\lambda }\left( z\right)
dz\leq \frac{1}{2\lambda }\dint\limits_{B}^{B+\rho }\left\vert q\left(
z\right) \right\vert \varphi _{\lambda }\left( z\right) dz.
\end{equation*}

To solve BVP3, we minimize the following cost functional: 
\begin{equation}
J_{\lambda }\left( W^{h}\right) =\left\Vert \left[ M_{N}W^{h}+P\left(
W_{x}^{h},W_{y}^{h},G_{x}^{h},G_{y}^{h},\mathbf{x}^{h}\right) \right]
e^{\lambda z}\right\Vert _{L_{2,N}^{h}\left( \Omega ^{h}\right) }^{2}.
\label{2.49}
\end{equation}

\textbf{Remark 2.2.} \emph{There are two differences between the functional }%
$J_{\lambda }\left( W^{h}\right) $\emph{\ and the analogous functional in 
\cite{kin1}, \cite[Chapter 11]{KL}. First, unlike these references, we do
not arrange zero boundary condition here for an analog of }$W^{h}$\emph{\
via sort of \textquotedblleft subtracting boundary conditions". Second, we
do not use a regularization penalization term in (\ref{2.49}). These two
differences simplify the current version of the convexification method, as
compared with the one in \cite{kin1}, \cite[Chapter 11]{KL}. Taken into
account the above mentioned equivalence of norms (\ref{2.40}) and (\ref%
{2.400}), the proof of the global strict convexity result, which is Theorem
3.1, is completely similar with the proof of Theorem 8.1 of \cite{kin1} and
Theorem 11.8.1 of \cite{KL}. Lemma 2.1 is used in this proof essentially.
Therefore, we omit below the proof of Theorem 3.1.}

To solve BVP3, we solve below the following minimization problem:

\textbf{Minimization Problem. }\emph{Minimize functional (\ref{2.49})} \emph{%
on the set }$\overline{K^{h}\left( R\right) }$\emph{\ defined in (\ref{2.47}%
).}

\section{Convergence Analysis}

\label{sec:3}

In this section, we formulate theorems 3.1-3.4 of the convergence analysis.
Keeping in mind (\ref{2.33}), we define the space $L_{2}^{h}\left( \Omega
_{1}^{h}\times \left( a,b\right) \right) $ similarly with (\ref{2.40}).

\textbf{Theorem 3.1 }(global strict convexity).

1. \emph{The functional }$J_{\lambda }\left( W^{h}\right) $\emph{\ has the Fr%
\'{e}chet derivative }$J_{\lambda }^{\prime }\left( W^{h}\right) $\emph{\ at
every point }$W^{h}\in K^{h}\left( 2R\right) $ \emph{and for all }$\lambda >0
$\emph{\ and} $J_{\lambda }^{\prime }\left( W^{h}\right) \in
H_{0,N}^{2,h}\left( \overline{\Omega }^{h}\right) .$ \emph{Hence, by (\ref%
{2.40})-(\ref{2.401})} $J_{\lambda }^{\prime }\left( W^{h}\right) \in
L_{2,N}^{h}\left( \Omega ^{h}\right) $\emph{\ as well. Furthermore, the Fr%
\'{e}chet derivative }$J_{\lambda }^{\prime }\left( W^{h}\right) $ \emph{%
satisfies Lipschitz condition on }$K^{h}\left( 2R\right) ,$\emph{\ i.e.
there exists a number }%
\begin{equation*}
\overline{C}=\overline{C}\left( h_{0},B,\rho ,\beta ,R,N,\Omega
^{h},\left\Vert g^{h}\mid _{D_{B+\rho }^{h}}\right\Vert _{L_{2}^{h}\left(
\Omega _{1}^{h}\times \left( a,b\right) \right) },\lambda \right) >0\emph{\ }
\end{equation*}%
\emph{depending only on listed parameters such that the following estimate
holds:}%
\begin{equation*}
\left\Vert J_{\lambda }^{\prime }\left( W_{2}^{h}\right) -J_{\lambda
}^{\prime }\left( W_{1}^{h}\right) \right\Vert _{L_{2,N}^{h}\left( \Omega
^{h}\right) }\leq \overline{C}\left\Vert W_{2}^{h}-W_{1}^{h}\right\Vert
_{H_{N}^{1,h}\left( \Omega ^{h}\right) },\text{ }\forall
W_{1}^{h},W_{2}^{h}\in K^{h}\left( 2R\right) .
\end{equation*}

2. \emph{There exist a sufficiently large number }$\lambda _{0}\geq 1$ \emph{%
and a number} $C>0,$%
\begin{equation}
\lambda _{0}=\lambda _{0}\left( h_{0},B,\rho R,N,c_{0},\Omega
^{h},\left\Vert g^{h}\mid _{D_{B+\rho }^{h}}\right\Vert _{L_{2}^{h}\left(
\Omega _{1}^{h}\times \left( a,b\right) \right) }\right) \geq 1,  \label{3.0}
\end{equation}%
\begin{equation}
C=C\left( h_{0},B,\rho ,R,N,c_{0},\Omega ^{h},\left\Vert g^{h}\mid
_{D_{B+\rho }^{h}}\right\Vert _{L_{2}^{h}\left( \Omega _{1}^{h}\times \left(
a,b\right) \right) }\right) >0,  \label{3.00}
\end{equation}%
\emph{\ \ both numbers depending} \emph{only on listed parameters, such that
for every }$\lambda \geq \lambda _{0}$\emph{\ the functional }$J_{\lambda
}\left( W^{h}\right) $\emph{\ is strictly convex on the closed set }$%
\overline{K^{h}\left( R\right) },$ \emph{i.e.} \emph{the following estimate
holds for all }$W_{1}^{h},W_{2}^{h}\in \overline{K^{h}\left( R\right) }:$%
\emph{\ \ }%
\begin{equation}
J_{\lambda }\left( W_{2}^{h}\right) -J_{\lambda }\left( W_{1}^{h}\right)
-J_{\lambda }^{\prime }\left( W_{1}^{h}\right) \left(
W_{2}^{h}-W_{1}^{h}\right) \geq C\left\Vert W_{2}^{h}-W_{1}^{h}\right\Vert
_{H_{2,N}^{1,h}\left( \Omega ^{h}\right) }^{2}.  \label{3.1}
\end{equation}

Below $C$ denotes different positive numbers depending only on parameters
listed in (\ref{3.00}).

\textbf{Theorem 3.2}.\emph{\ Let }$\lambda _{0}\geq 1$\emph{\ be the number
of Theorem 3.1. Then for every }$\lambda \geq \lambda _{0}$\emph{\ there
exists a single minimizer }$W_{\lambda ,\min }^{h}\in \overline{K^{h}\left(
R\right) }$\emph{\ of the functional }$J_{\lambda }\left( W^{h}\right) $%
\emph{\ on the set }$\overline{K^{h}\left( R\right) }$\emph{\ and the
following inequality holds:} 
\begin{equation}
-J_{\lambda }^{\prime }\left( W_{\lambda ,\min }^{h}\right) \left(
W^{h}-W_{\lambda ,\min }^{h}\right) \leq 0,\text{ }\forall W^{h}\in 
\overline{K^{h}\left( R\right) }.  \label{3.2}
\end{equation}

According to the regularization theory \cite{BK,T}, the minimizer $%
W_{\lambda ,\min }^{h}$ of functional (\ref{2.49}) is called
\textquotedblleft regularized solution". It is important to estimate the
accuracy of the regularized solution depending on the noise in the data. To
do this, we recall first that, following the regularization theory, we need
to assume the existence of the \textquotedblleft ideal" solution of BVP3,
i.e. solution with the noiseless data. The ideal solution is also called
\textquotedblleft exact" solution. We denote this solution $W^{h\ast }\in
H_{N}^{2,h}\left( \Omega ^{h}\right) .$ We denote the noiseless data in (\ref%
{2.45}), (\ref{2.46}) as $\widetilde{G}^{h\ast },g^{h\ast }.$ We assume that 
\begin{equation}
W^{h\ast }\in K^{h\ast }\left( R\right) =\left\{ 
\begin{array}{c}
W^{h}\in H_{N}^{2,h}\left( \Omega ^{h}\right) :u^{h}\left( \mathbf{x},\alpha
\right) \geq c_{0}^{2}\text{ in }\overline{\Omega }^{h}, \\ 
W^{h}\mid _{\Gamma ^{h}\cup D_{B}^{h}}=\widetilde{G}^{h\ast }\left( \mathbf{x%
}\right) ,\left\Vert W^{h}\right\Vert _{H_{N}^{1,h}\left( \Omega ^{h}\right)
}<R%
\end{array}%
\right\} ,  \label{3.3}
\end{equation}%
see (\ref{2.47}). Let $m^{h\ast }\left( \mathbf{x}^{h}\right) $ be the
semidiscrete exact target function which is found via the substitution of $%
W^{h,\ast }$ first in (\ref{2.42}) and (\ref{2.43}) and then in (\ref{2.160}%
).

In the reality, however, the data $g^{h}$ and $\widetilde{G}^{h}$ are always
noisy. Let $\delta \in \left( 0,1\right) $ be the level of the noise in the
data. We assume that 
\begin{equation*}
\left\Vert \left( g^{h}-g^{h\ast }\right) \left( \mathbf{x}^{h},\alpha
\right) \mid _{D_{B+\rho }^{h}}\right\Vert _{L_{2}^{h}\left( \Omega
_{1}^{h}\times \left( a,b\right) \right) }+
\end{equation*}%
\begin{equation}
+\left\Vert \partial _{\alpha }\left( g^{h}-g^{h\ast }\right) \left( \mathbf{%
x}^{h},\alpha \right) \mid _{D_{B+\rho }^{h}}\right\Vert _{L_{2}^{h}\left(
\Omega _{1}^{h}\times \left( a,b\right) \right) }<\delta .  \label{3.4}
\end{equation}%
Let $F^{h\ast }\left( \mathbf{x}^{h}\right) \in H_{N}^{1,h}\left( \Omega
^{h}\right) $ be an extension of the boundary vector function $\widetilde{G}%
^{h\ast }$ from $\Gamma ^{h}\cup D_{B}^{h}$ inside of the semidiscrete
domain $\Omega ^{h}.$ Such an extension exists since the ideal solution $%
W^{h\ast }\in H_{N}^{1,h}\left( \Omega ^{h}\right) $ and $\widetilde{G}%
^{h\ast }$ is its boundary condition. Thus, $F^{h\ast }\mid _{\Gamma
^{h}\cup D_{B}^{h}}=\widetilde{G}^{h\ast }.$ We assume that there exists an
extension $F^{h}\in H_{N}^{1,h}\left( \Omega ^{h}\right) $ of the noisy
boundary vector function $\widetilde{G}^{h}$ from the boundary $\Gamma
^{h}\cup D_{B}^{h}$ inside the domain $\Omega ^{h}$ such that $F^{h}\mid
_{\Gamma ^{h}\cup D_{B}^{h}}=\widetilde{G}^{h}$ and 
\begin{equation}
\left\Vert F^{h}-F^{h\ast }\right\Vert _{H_{N}^{1,h}\left( \Omega
^{h}\right) }<\delta .  \label{3.5}
\end{equation}%
In addition, we assume that 
\begin{equation}
\left\Vert F^{h}\right\Vert _{H_{N}^{1,h}\left( \Omega ^{h}\right)
},\left\Vert F^{h\ast }\right\Vert _{H_{N}^{1,h}\left( \Omega ^{h}\right)
}<R.  \label{3.6}
\end{equation}

Theorem 3.3 provides the desired accuracy estimate of the regularized
solution $W_{\lambda ,\min }^{h}$ depending on the level $\delta $ of the
noise in the data.

\textbf{Theorem 3.3. }\emph{Assume that conditions (\ref{3.3})-(\ref{3.6})
hold. Furthermore, assume that }%
\begin{equation}
\left\Vert W^{h\ast }\right\Vert _{H_{N}^{1,h}\left( \Omega ^{h}\right)
}<R-\alpha ,  \label{3.7}
\end{equation}%
\emph{where the number }$\alpha \in \left( 0,R\right) $\emph{\ is so small
that }%
\begin{equation}
\alpha <C\delta .  \label{3.8}
\end{equation}%
\emph{Let }$\lambda _{1}$\emph{\ be the number }$\lambda _{0}$\emph{\ of
Theorem 3.1 in the case when }$R$\emph{\ in (\ref{3.0}) is replaced with }$%
2R,$%
\begin{equation}
\lambda _{1}=\lambda _{0}\left( h_{0},B,\rho ,2R,N,c_{0},\Omega
^{h},\left\Vert g^{h}\mid _{D_{B+\rho }^{h}}\right\Vert _{L_{2}^{h}\left(
\Omega _{1}^{h}\times \left( a,b\right) \right) }\right) \geq 1.
\label{3.80}
\end{equation}%
\emph{Let }$W_{\lambda _{1},\min }^{h}\in \overline{K^{h}\left( R\right) }$%
\emph{\ be the minimizer on of the functional }$J_{\lambda }\left(
W^{h}\right) $\emph{\ on the set }$\overline{K^{h}\left( R\right) },$\emph{\
which is claimed in Theorem 3.2. Then the following accuracy estimate holds:}%
\begin{equation}
\left\Vert W_{\lambda _{1},\min }^{h}-W^{h\ast }\right\Vert
_{H_{N}^{1,h}\left( \Omega ^{h}\right) }\leq C\delta .  \label{3.9}
\end{equation}

\textbf{Remark 3.1.} \emph{It easily follows from the proof of this theorem
that the estimate similar with the one in (\ref{3.9}) is valid for any }$%
\lambda \geq \lambda _{1}.$\emph{\ However, in this case one needs to
replace in (\ref{3.9}) }$W_{\lambda _{1},\min }^{h}$\emph{\ with }$%
W_{\lambda ,\min }^{h}$\emph{\ and the right hand side should be replaced
with }$Ce^{2\lambda \left( B+\rho \right) }\delta .$\emph{\ A similar
statement is true for the gradient descent method formulated below.}

We now construct the gradient descent method of the minimization of the
functional $J_{\lambda _{1}}\left( W^{h}\right) .$ Let $W_{0}^{h}\in
K^{h}\left( R/3\right) $ be an arbitrary point and $\gamma \in \left(
0,1\right) $ be a number. The sequence of the gradient descent method is:%
\begin{equation}
W_{n}^{h}=W_{n-1}^{h}-\gamma J_{\lambda _{1}}\left( W_{n-1}^{h}\right) ,%
\text{ }n=1,2,...  \label{3.10}
\end{equation}%
Note that since by Theorem 3.1 $J_{\lambda _{1}}^{\prime }\left(
W_{n-1}^{h}\right) \in H_{0,N}^{1,h}\left( \overline{\Omega }^{h}\right) ,$ $%
\forall n$, then all vector functions $W_{n}^{h}$ satisfy the same boundary
condition as the one in (\ref{2.47}).

\textbf{Theorem 3.4. }\emph{Let }$C\delta \in \left( 0,R/3\right) $\emph{\
and let the number }$\beta \in \left( C\delta ,R/3\right) .$\emph{\ Suppose
that the exact solution }$W^{h\ast }\in K^{h}\left( R/3-\beta \right) $\emph{%
. Let }$\lambda =\lambda _{1}$\emph{\ where }$\lambda _{1}$ \emph{is defined
in (\ref{3.80}). Then there exists a sufficiently small number }$\gamma
_{0}\in \left( 0,1\right) $\emph{\ such that for any }$\gamma \in \left(
0,\gamma _{0}\right) $\emph{\ all terms of the sequence (\ref{3.10}) }$%
W_{n}^{h}\in K^{h}\left( R/3\right) .$\emph{\ Furthermore, there exists a
number }$\theta =\theta \left( \gamma \right) \in \left( 0,\gamma \right) $%
\emph{\ such that the following convergence estimate holds}%
\begin{equation}
\left\Vert W_{n}^{h}-W^{h\ast }\right\Vert _{H_{N}^{1,h}\left( \Omega
^{h}\right) }\leq C\delta +\theta ^{n}\left\Vert W_{0}^{h}-W_{\lambda
_{1},\min }^{h}\right\Vert _{H_{N}^{1,h}\left( \Omega ^{h}\right) }.
\label{3.11}
\end{equation}%
\emph{In addition to the function }$m^{h\ast }\left( \mathbf{x}^{h}\right) ,$%
\emph{\ let }$m_{n}^{h}\left( \mathbf{x}^{h}\right) $\emph{\ be the
semidiscrete target function, which is found via the substitution of }$%
W_{n}^{h}$ \emph{\ first in (\ref{2.42}) and (\ref{2.43}) and then in the
left hand side of (\ref{2.160}). Then }%
\begin{equation}
\left\Vert m_{n}^{h}-m^{h\ast }\right\Vert _{H_{N}^{1,h}\left( \Omega
^{h}\right) }\leq C\delta +\theta ^{n}\left\Vert W_{0}^{h}-W_{\lambda
_{1},\min }^{h}\right\Vert _{H_{N}^{1,h}\left( \Omega ^{h}\right) }.
\label{3.110}
\end{equation}

\textbf{Remark 3.2.}\emph{\ Since a smallness assumption is not imposed on
the number }$R$\emph{\ and since the starting point }$W_{0}^{h}$\emph{\ of
the gradient descent method (\ref{3.10}) is an arbitrary point of the set }$%
K^{h}\left( R\right) ,$\emph{\ then Theorem 3.4 claims the global
convergence of our method, see section 1 for our definition of the global
convergence. }

\section{Proofs}

\label{sec:4}

We omit the proof of Theorem 3.1, see Remark 2.2. Theorem 3.2 follows
immediately from Theorem 3.1 and a combination of Lemma 2.1 and Theorem 2.1
of \cite{Bak}. Hence, we prove only Theorems 3.3 and 3.4.

\subsection{Proof of Theorem 3.3}

\label{sec:4.1}

Denote%
\begin{equation*}
\widetilde{K}_{0}^{h}\left( 2R\right) =\left\{ V^{h}\in H_{0,N}^{1,h}\left(
\Omega ^{h}\right) :\left\Vert V^{h}\right\Vert _{H_{N}^{1,h}\left( \Omega
^{h}\right) }<2R\right\} ,
\end{equation*}%
\begin{equation}
V^{h}=W^{h}-F^{h},\text{ }\forall W^{h}\in K^{h}\left( R\right) ,
\label{4.1}
\end{equation}%
\begin{equation}
V^{h\ast }=W^{h\ast }-F^{h\ast }.  \label{4.2}
\end{equation}%
By (\ref{3.6}), (\ref{4.1}) and (\ref{4.2}) 
\begin{equation}
\left\Vert V^{h\ast }\right\Vert _{H_{N}^{1,h}\left( \Omega ^{h}\right)
}<2R,\left\Vert V^{h}\right\Vert _{H_{N}^{1,h}\left( \Omega ^{h}\right) }<2R.
\label{4.3}
\end{equation}%
For each vector function $V^{h}\in $ $\widetilde{K}_{0}^{h}\left( 2R\right) ,
$ consider the vector function $W^{h}\left( V^{h}\right) =V^{h}+F^{h}$ and
then, using $W^{h}\left( V^{h}\right) ,$ construct the function $u\left(
V^{h}\right) \left( \mathbf{x},\alpha \right) $ as in (\ref{2.42}), (\ref%
{2.43}). Hence, we denote%
\begin{equation}
K_{0}^{h}\left( 2R\right) =\left\{ 
\begin{array}{c}
V^{h}\in H_{0,N}^{1,h}\left( \Omega ^{h}\right) : \\ 
u\left( V^{h}+F^{h}\right) \left( \mathbf{x},\alpha \right) \geq c_{0}^{2}%
\text{ in }\overline{\Omega }^{h},\left\Vert V^{h}\right\Vert
_{H_{N}^{1,h}\left( \Omega ^{h}\right) }<2R%
\end{array}%
\right\} .  \label{4.4}
\end{equation}%
By (\ref{2.47}), (\ref{3.3}), (\ref{4.3}) and (\ref{4.4})%
\begin{equation}
V^{h}=\left( W^{h}-F^{h}\right) \in K_{0}^{h}\left( 2R\right) ,\text{ }%
\forall W^{h}\in K^{h}\left( R\right) ,  \label{4.5}
\end{equation}%
\begin{equation}
V^{h\ast }=\left( W^{h\ast }-F^{h\ast }\right) \in K_{0}^{h}\left( 2R\right)
.  \label{4.6}
\end{equation}%
Consider the functional $I_{\lambda _{1}}:K_{0}^{h}\left( 2R\right)
\rightarrow \mathbb{R},$ where 
\begin{equation*}
\text{ }I_{\lambda _{1}}\left( V^{h}\right) =J_{\lambda _{1}}\left(
V^{h}+F^{h}\right) .
\end{equation*}%
Then obvious analogs of Theorems 3.1, 3.2 are valid for $I_{\lambda
_{1}}\left( V^{h}\right) $ with the replacement of the pair $\left(
K^{h}\left( R\right) ,\lambda _{0}\right) $ with with the pair $\left(
K_{0}^{h}\left( 2R\right) ,\lambda _{1}\right) ,$ where $\lambda _{1}$ is
defined in (\ref{3.80}). Let $V_{\lambda _{1},\min }^{h}\in \overline{%
K_{0}^{h}\left( 2R\right) }$ be the minimizer of the functional $I_{\lambda
_{1}}\left( V^{h}\right) $ on the set $\overline{K_{0}^{h}\left( 2R\right) }$%
. The existence and uniqueness of this minimizer follows from that analog of
Theorem 3.2. Since by (\ref{4.5}) and (\ref{4.6}) both vector functions $%
V^{h},V^{h\ast }\in K_{0}^{h}\left( 2R\right) ,$ then (\ref{3.1}) implies 
\begin{equation}
J_{\lambda _{1}}\left( V^{h\ast }+F^{h}\right) -J_{\lambda _{1}}\left(
V_{\lambda _{1},\min }^{h}+F^{h}\right) -J_{\lambda _{1}}^{\prime }\left(
V_{\lambda _{1},\min }^{h}+F^{h}\right) \left( V^{h\ast }-V_{\lambda
_{1},\min }^{h}\right) \geq   \label{4.7}
\end{equation}%
\begin{equation*}
\geq C\left\Vert V_{\lambda _{1},\min }^{h}-V^{h\ast }\right\Vert
_{H_{N}^{1,h}\left( \Omega ^{h}\right) }^{2}.
\end{equation*}%
By (\ref{3.2})%
\begin{equation*}
-J_{\lambda _{1}}^{\prime }\left( V_{\lambda _{1},\min }^{h}+F^{h}\right)
\left( V^{h\ast }-V_{\lambda _{1},\min }^{h}\right) \leq 0.
\end{equation*}%
Hence, (\ref{4.7}) implies%
\begin{equation}
J_{\lambda _{1}}\left( V^{h\ast }+F^{h}\right) \geq C\left\Vert V_{\lambda
_{1},\min }^{h}-V^{h\ast }\right\Vert _{H_{N}^{1,h}\left( \Omega ^{h}\right)
}^{2}.  \label{4.8}
\end{equation}%
Estimate now the left hand side of (\ref{4.8}) from the above. Since the
vector function $W^{h\ast }=V^{h\ast }+F^{h\ast }$ satisfies equation (\ref%
{2.44}), then (\ref{2.49}) implies that $J_{\lambda _{1}}\left( V^{h\ast
}+F^{h\ast }\right) =J_{\lambda _{1}}\left( W^{h\ast }\right) =0.$ Hence,
using (\ref{2.49}) and (\ref{3.5}), we obtain 
\begin{equation*}
J_{\lambda _{1}}\left( V^{h\ast }+F^{h}\right) =J_{\lambda _{1}}\left(
W^{h\ast }+\left( F^{h}-F^{h\ast }\right) \right) \leq Ce^{2\lambda
_{1}\left( B+\rho \right) }\delta ^{2}.
\end{equation*}%
Hence, (\ref{4.8}) implies that with a different constant $C$%
\begin{equation}
\left\Vert V_{\lambda _{1},\min }^{h}-V^{h\ast }\right\Vert
_{H_{N}^{1,h}\left( \Omega ^{h}\right) }\leq C\delta .  \label{4.9}
\end{equation}%
Let 
\begin{equation}
\overline{W}_{\lambda _{1},\min }^{h}=V_{\lambda _{1},\min }^{h}+F^{h}.
\label{4.10}
\end{equation}%
Using (\ref{4.9}), (\ref{4.10}) and the triangle inequality, we obtain 
\begin{equation}
\left\Vert \overline{W}_{\lambda _{1},\min }^{h}-W^{h\ast }\right\Vert
_{H_{N}^{1,h}\left( \Omega ^{h}\right) }\leq C\delta .  \label{4.100}
\end{equation}%
Hence, (\ref{3.7}), (\ref{3.8}) and (\ref{4.100}) imply 
\begin{equation}
\overline{W}_{\lambda _{1},\min }^{h}\in K^{h}\left( R\right) .  \label{4.11}
\end{equation}

Consider now the minimizer $W_{\lambda _{1},\min }^{h}\in \overline{%
K^{h}\left( R\right) }$ of the functional $J_{\lambda _{1}}\left(
W^{h}\right) $ on the set $\overline{K^{h}\left( R\right) },$ which is
claimed by Theorem 3.1. Then by (\ref{4.11}) 
\begin{equation}
J_{\lambda _{1}}\left( W_{\lambda _{1},\min }^{h}\right) \leq J_{\lambda
_{1}}\left( \overline{W}_{\lambda _{1},\min }^{h}\right) .  \label{4.12}
\end{equation}%
On the other hand, let $\widetilde{V}_{\lambda _{1},\min }^{h}=W_{\lambda
_{1},\min }^{h}-F^{h}.$ Then $\widetilde{V}_{\lambda _{1},\min }^{h}\in
K_{0}^{h}\left( 2R\right) .$ Hence, 
\begin{equation*}
J_{\lambda _{1}}\left( W_{\lambda _{1},\min }^{h}\right) =J_{\lambda
_{1}}\left( \widetilde{V}_{\lambda _{1},\min }^{h}+F^{h}\right) \geq
J_{\lambda _{1}}\left( V_{\lambda _{1},\min }^{h}+F^{h}\right) =J_{\lambda
_{1}}\left( \overline{W}_{\lambda _{1},\min }^{h}\right) .
\end{equation*}%
Thus, $J_{\lambda _{1}}\left( W_{\lambda _{1},\min }^{h}\right) \geq
J_{\lambda _{1}}\left( \overline{W}_{\lambda _{1},\min }^{h}\right) .$
Comparing this with (\ref{4.12}), we obtain $J_{\lambda _{1}}\left(
W_{\lambda _{1},\min }^{h}\right) =J_{\lambda _{1}}\left( \overline{W}%
_{\lambda _{1},\min }^{h}\right) .$ Since by Theorem 3.2 the minimizer of\
the functional $J_{\lambda _{1}}\left( W^{h}\right) $ on the set $\overline{%
K^{h}\left( R\right) }$ is unique, then $W_{\lambda _{1},\min }^{h}=%
\overline{W}_{\lambda _{1},\min }^{h}.$ This and (\ref{4.100}) prove (\ref%
{3.9}). $\square $

\subsection{Proof of Theorem 3.4}

\label{sec:4.2}

\bigskip By the triangle inequality and (\ref{3.9}) 
\begin{equation}
\left\Vert W_{\lambda _{1},\min }^{h}\right\Vert _{H_{N}^{1,h}\left( \Omega
^{h}\right) }-\left\Vert W^{h\ast }\right\Vert _{H_{N}^{1,h}\left( \Omega
^{h}\right) }\leq \left\Vert W_{\lambda _{1},\min }^{h}-W^{h\ast
}\right\Vert _{H_{N}^{1,h}\left( \Omega ^{h}\right) }\leq C\delta .
\label{4.13}
\end{equation}%
By (\ref{4.13}) 
\begin{equation}
\left\Vert W_{\lambda _{1},\min }^{h}\right\Vert _{H_{N}^{1,h}\left( \Omega
^{h}\right) }\leq \left\Vert W^{h\ast }\right\Vert _{H_{N}^{1,h}\left(
\Omega ^{h}\right) }+C\delta \leq \frac{R}{3}-\left( \beta -C\delta \right) <%
\frac{R}{3}.  \label{4.14}
\end{equation}%
Since the starting point of sequence (\ref{3.10}) $W_{0}^{h}\in K^{h}\left(
R/3\right) ,$ then (\ref{4.14}) and Theorem 6 of \cite{SAR2} imply that
there exists a sufficiently small number $\gamma _{0}\in \left( 0,1\right) $
such that for every $\gamma \in \left( 0,\gamma _{0}\right) $ all vector
functions $W_{n}^{h}\in K^{h}\left( R/3\right) ,n=1,2,...$ and also that
there exists a number $\theta =\theta \left( \gamma \right) \in \left(
0,1\right) $ such that 
\begin{equation}
\left\Vert W_{n}^{h}-W_{\lambda _{1},\min }^{h}\right\Vert
_{H_{N}^{1,h}\left( \Omega ^{h}\right) }\leq \theta ^{n}\left\Vert
W_{0}^{h}-W_{\lambda _{1},\min }^{h}\right\Vert _{H_{N}^{1,h}\left( \Omega
^{h}\right) },\text{ }n=1,2,...  \label{4.15}
\end{equation}%
Next, using (\ref{3.9}), (\ref{4.15}) and the triangle inequality, we obtain 
\begin{equation*}
\left\Vert W_{n}^{h}-W^{h\ast }\right\Vert _{H_{N}^{1,h}\left( \Omega
^{h}\right) }\leq \left\Vert W_{\lambda _{1},\min }^{h}-W^{h\ast
}\right\Vert _{H_{N}^{1,h}\left( \Omega ^{h}\right) }+\left\Vert
W_{n}^{h}-W_{\lambda _{1},\min }^{h}\right\Vert _{H_{N}^{1,h}\left( \Omega
^{h}\right) }\leq
\end{equation*}%
\begin{equation*}
\leq C\delta +\theta ^{n}\left\Vert W_{0}^{h}-W_{\lambda _{1},\min
}^{h}\right\Vert _{H_{N}^{1,h}\left( \Omega ^{h}\right) },\text{ }n=1,2,...,
\end{equation*}%
which proves (\ref{3.11}). The final convergence estimate (\ref{3.110}) of
this theorem follows immediately from (\ref{3.11}) and the above described
procedure of the construction of functions $m^{\ast }\left( \mathbf{x}%
^{h}\right) $ and $m_{n}\left( \mathbf{x}^{h}\right) $ from vector functions 
$W^{h\ast }\left( \mathbf{x}^{h}\right) $ and $W_{n}^{h}\left( \mathbf{x}%
^{h}\right) $ respectively. $\square $

\section{Numerical studies}

\label{5}

\subsection{Numerical implementation}

\label{5.1}

\begin{figure}[tbp]
\begin{center}
\begin{tabular}{c}
\includegraphics[width=5cm]{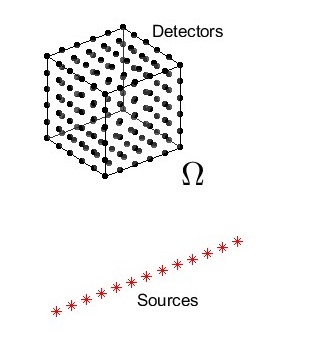}%
\end{tabular}%
\end{center}
\caption{\emph{A schematic diagram
measurements.}}
\label{example0}
\end{figure}

In all the numerical tests, we have chosen in {(\ref{2.1}) }the numbers $B=0$%
, $\rho =1$, i.e. the domain $\Omega =[0,1]^{3}$ in {(\ref{2.1})} is a unit
cube. We set in {(\ref{2.8}) }$a=-2$, $b=3$, $d=1/2$, $z_{0}=-1$. We have
used 101 sources, which were uniformly distributed on the line $L_{src}$ and
detectors were uniformly distributed on the surface $\partial \Omega $ with
mesh size $1/20\times 1/20$. See Figure \ref{example0} for a schematic
diagram of our measurements. 

\subsubsection{The forward problem of TTTP (\protect\ref{1.3})}

\label{5.1.1}

We solve the forward problem to generate the data of the arriving time $%
g\left( \mathbf{x},\mathbf{x}_{0}\right) =\tau \left( \mathbf{x},\mathbf{x}%
_{0}\right) $ in (\ref{1.3}) for TTTP of section 1 as follows. Given the
heterogeneous medium $m(\mathbf{x})$, instead of solving the nonlinear
eikonal equation \eqref{1.2} and \eqref{1.20} directly, we have solved the
forward problem (\ref{1.2})-(\ref{1.20}) for each of those 101 sources
located on $L_{src}$ by the 3D fast marching method built in Matlab.  The
Fast Marching method is a technique that approximates the solution of the Eikonal nonlinear partial differential equation, it is similar to the Dijkstra algorithm to find the shortest paths in graphs which are uniformly sized spatial grid with sound speed value at each node.  For each source located on $L_{src}$, the Fast
Marching method simultaneously computes the shortest paths $\Gamma (\mathbf{x%
},\mathbf{x}_{0})$ on graphs and the arriving time $\tau (\mathbf{x},\mathbf{%
x}_{0})$. Here we set the step size with respect to spatial variables for
the fast marching method $h_{\text{sp}}=1/30$.

\subsubsection{The inverse problem}

\label{5.1.2} Having the computationally simulated data {(\ref{1.3}) for the
inverse problem, which are computed }as in sub-subsection \ref{5.1.1}, we
add the random noise to these data as: 
\begin{equation*}
g_{\text{noise}}\left( \mathbf{x},\mathbf{x}_{0}\right) =g\left( \mathbf{x},%
\mathbf{x}_{0}\right) +\delta \max_{\mathbf{x}\in \partial \Omega ,\mathbf{x}%
_{0}\in L_{\text{src}}}|g\left( \mathbf{x},\mathbf{x}_{0}\right) |\xi _{%
\mathbf{x}_{0}},\text{ }\mathbf{x}\in \partial \Omega ,\mathbf{x}_{0}\in L_{%
\text{src}},
\end{equation*}%
Here $\xi _{\mathbf{x}_{0}}\in \left( -1,1\right) $ is the uniformly
distributed random variable and $\delta =0.05$, i.e. $5\%$ noise level in
all the numerical tests.

{{To solve inverse problem \eqref{1.3} numerically, we have first selected
optimal values for some parameters. Those values were selected by the trial
and error procedure. We show below the effects of different combinations of
these parameters in the first two tests. As soon as the best parameters were
selected in the first two tests, they were used in the rest of tests. }}

{We have solved the following minimization problem:}

\textbf{Minimization Problem.}{\ \emph{Minimize the functional} } 
\begin{equation}
J_{\lambda ,\beta }\left( W^{h}\right) =\left\Vert \left[ M_{N}W^{h}+P\left(
W_{x}^{h},W_{y}^{h},g_{x}^{h},g_{y}^{h},\mathbf{x}^{h}\right) \right]
e^{\lambda z}\right\Vert _{L_{2,N}^{h}\left( \Omega ^{h}\right) }^{2}+\beta
\left\Vert W^{h}\right\Vert _{H_{N}^{2,h}\left( \Omega ^{h}\right) }^{2},
\label{mini}
\end{equation}%
\emph{on the set }$\overline{K^{h}\left( R\right) }$\emph{\ defined in (\ref%
{2.47}), where }$\beta \in \left( 0,1\right) $\emph{\ is the regularization
parameter.}

Here the space $H_{N}^{2,h}\left( \Omega ^{h}\right) $ is defined as%
\begin{equation*}
H_{N}^{2,h}\left( \Omega ^{h}\right) =\left\{ 
\begin{array}{c}
Q^{h}(\mathbf{x}^{h}):\left\Vert Q^{h}\right\Vert _{H_{N}^{2,h}\left( \Omega
^{h}\right) }^{2}= \\ 
\left\Vert Q_{x}^{h}\right\Vert _{L_{2N}^{h}\left( \Omega ^{h}\right)
}^{2}+\left\Vert Q_{y}^{h}\right\Vert _{L_{2N}^{h}\left( \Omega ^{h}\right)
}^{2}+\dsum\limits_{s=0}^{2}\left\Vert \partial _{z}^{s}Q^{h}\right\Vert
_{L_{2N}^{h}\left( \Omega ^{h}\right) }^{2}<\infty%
\end{array}%
\right\} .
\end{equation*}

The regularization term $\beta \left\Vert W^{h}\right\Vert
_{H_{N}^{2,h}\left( \Omega ^{h}\right) }^{2}$ in {{(\ref{mini})}} is not
involved in our above theory. On the other hand, we see in our numerical
experiments that our method does not perform well without the regularization
term. We cannot explain yet why this takes place, and this should be a
subject of our further research. We note that since the regularization term
represents a strictly convex functional $Y_{\beta }:H_{N}^{2,h}\left( \Omega
^{h}\right) \rightarrow \mathbb{R},$ 
\begin{equation*}
Y_{\beta }\left( W^{h}\right) =\beta \left\Vert W^{h}\right\Vert
_{H_{N}^{2,h}\left( \Omega ^{h}\right) }^{2},
\end{equation*}%
then an obvious analog of Theorem 3.1 implies that the functional $%
J_{\lambda ,\beta }\left( W^{h}\right) :K_{2}^{h}\left( R\right) \rightarrow 
\mathbb{R}$ is strictly convex on the set $\overline{K_{2}^{h}\left(
R\right) },$ where $K_{2}^{h}\left( R\right) $ is obtained from the set $%
K^{h}\left( R\right) $ in {{(\ref{2.47}) via replacing in (\ref{2.47}) }}$%
H_{N}^{1,h}\left( \Omega ^{h}\right) $ with $H_{N}^{2,h}\left( \Omega
^{h}\right) .$ Then {{(\ref{3.1}) }}should be replaced with%
\begin{equation*}
J_{\lambda ,\beta }\left( W_{2}^{h}\right) -J_{\lambda ,\beta }\left(
W_{1}^{h}\right) -J_{\lambda ,\beta }^{\prime }\left( W_{1}^{h}\right)
\left( W_{2}^{h}-W_{1}^{h}\right) \geq 
\end{equation*}%
\begin{equation*}
\geq C\left\Vert W_{2}^{h}-W_{1}^{h}\right\Vert _{H_{2,N}^{1,h}\left( \Omega
^{h}\right) }^{2}+\beta \left\Vert W_{2}^{h}-W_{1}^{h}\right\Vert
_{H_{2,N}^{2,h}\left( \Omega ^{h}\right) }^{2},\forall
W_{1}^{h},W_{2}^{h}\in \overline{K_{2}^{h}\left( R\right) }.
\end{equation*}

We have minimized the fully discrete version of $J_{\lambda ,\beta }\left(
W^{h}\right) $ with respect to the values of the corresponding vector
function at grid points. To minimize the discretized functional $J_{\lambda
,\beta }\left( W^{h}\right) $, we use the Matlab's built-in function \textbf{%
fminunc} to solve the optimization problem. This function calculates the
gradient $\nabla J_{\lambda ,\beta }\left( W^{h}\right) $ automatically, and
we let the iterations stop when the condition $|\nabla J_{\lambda ,\beta
}\left( W^{h}\right) |<10^{-6}$ is fulfilled.

\subsection{Results}

\label{5.2}

In the numerical tests of this subsection, we demonstrate the efficiency of
our method. We display on Figures the true and computed functions $n\left( 
\mathbf{x}\right) =\sqrt{m\left( \mathbf{x}\right) }$, see section 1. In all
tests the background value $n_{\text{bkgr}}\left( \mathbf{x}\right) =1$,
which corresponds to the background value $c_{\text{bkgr}}\left( \mathbf{x}%
\right) =1.$ In those figures below, we depict 2-D slices to demonstrate the
values of the true function $n\left( \mathbf{x}\right) ${\ and computed
function }$n_{\text{comp}}\left( \mathbf{x}\right) ${. \ }

We have chosen values of all our parameters by the trial and error
procedure. In all tests, we set the regularization parameter in {{(\ref{mini}%
)}} $\beta =10^{-4}.$ In the first two tests, we use the mesh step size $%
h=1/10$. And for the tests number 3-5, we use the mesh step size $h=1/15$
for better resolutions. In Tests 1 and 2 we select an optimal pair $\left(
\lambda ,N\right) $ of parameters. Therefore, we work in the first two tests
with a ball-shaped inclusion which is a rather simple shape. As soon as an
optimal pair $\left( \lambda ,N\right) $ is selected, we work in Tests 3-5
with letter-like shapes of inclusions. We have chosen letters since they are
non convex and have voids, i.e. their shapes are rather complicated ones.

\textbf{Test 1}. First, we test the reconstruction by our method of the case
of a ball-shaped inclusion. The true function $n(\mathbf{x})$ is depicted on
Figures \ref{example1} (a) and (b), $n(\mathbf{x})=1.5$ inside of this
inclusion and $n(\mathbf{x})=1$ outside of it. Since $m(\mathbf{x})=1.5^{2}$
inside of this inclusion and $m(\mathbf{x})=1$ outside of it, then the
inclusion/background contrast in the target function $m(\mathbf{x})$ is $%
1.5^{2}:1=2.25:1.$ In this test, we fix $N=6$ and vary $\lambda $ from 0 to
4. See Figures \ref{example1} for the reconstruction results. One can
observe that the result with $\lambda =0$ i.e. in the case when the Carleman
Weight Function is absent in the functional $J_{\lambda ,\rho }\left(
W^{h}\right) ,$ is unacceptable. Even though results with $\lambda =2,3,4$
are about the same$,$ we choose $\lambda =4$ as the optimal value since our
theory basically says that larger values of $\lambda $ are better than lower
ones.

\textbf{Test 2}. In this case, we test the influence of the parameter $N$.
We again use the same inclusion as the one in Test 1. We fix the parameter $%
\lambda =4,$ which we have chosen in Test 1, and allow $N$ to vary from 4 to
10. See Figures \ref{example2} for the results of the reconstruction. One
can observe that the reconstruction results for $N=6,8,10$ are basically the
same. Therefore, we choose $N=6$ since this choice ensures a lesser
computational cost.

In conclusion, we select an optimal pair $\left( \lambda ,N\right) $ of
parameters as: 
\begin{equation}
\lambda =4,N=6\text{ in tests }3-5.  \label{5.10}
\end{equation}

In tests 3-5 inclusions are letter-shaped. We have intentionally chosen
these shapes since they are non convex and are, therefore, hard to image.

\textbf{Test 3}. We test the reconstruction by our method of the case when
the shape of our inclusion is the same as the shape of the letter `$A$'. The
function $n(\mathbf{x})$ is depicted on Figures \ref{example3} (a) and (b). $%
n=1.5$ inside of this inclusion and $n=1$ outside of it. See Figures \ref%
{example3} for the reconstruction results.

\textbf{Test 4}. We test the reconstruction by our method for the case when
the shape of our inclusion is the same as the shape of the letter `$C$'. The
function $n(\mathbf{x})$ is depicted on Figures \ref{example4} (a) and (b). $%
n=1.5$ inside of this inclusion and $n=1$ outside of it. See Figures \ref%
{example4} for the reconstruction results.

\textbf{Test 5}. We test the reconstruction by our method of the case when
the shape of our inclusion is the same as the shape of the letter `$\Omega $%
'. The function $n(\mathbf{x})$ is depicted on Figures \ref{example5} (a)
and (b). $n=1.5$ inside of this inclusion and $n=1$ outside of it. See
Figures \ref{example5} for the reconstruction results.

\section{Summary}

\label{6}

We have presented the first computational result for the Travel Time
Tomography Problem in the 3-D case. To do this, we have implemented
numerically the version of \cite{kin1}, \cite[Chapter 11]{KL} of the
globally convergent convexification numerical method. This method minimizes
a certain weighted cost functional with the Carleman Weight Function in it.
We have provided the global convergence analysis for our method.
Interestingly our computations for Test 1 with $\lambda =0$ show that
results have an unacceptable quality when the Carleman Weight Function is
absent in this functional. We have imaged letter-shaped inclusions, which
are non convex and are, therefore, hard to image. Nevertheless, shapes of
inclusions are imaged accurately in all tests.

\begin{figure}[tbp]
\begin{center}
\begin{tabular}{cc}
\includegraphics[width=3.5cm]{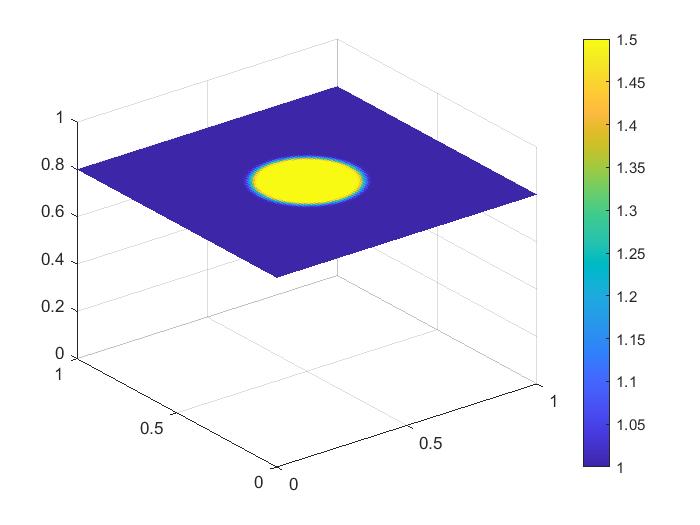} & %
\includegraphics[width=3.5cm]{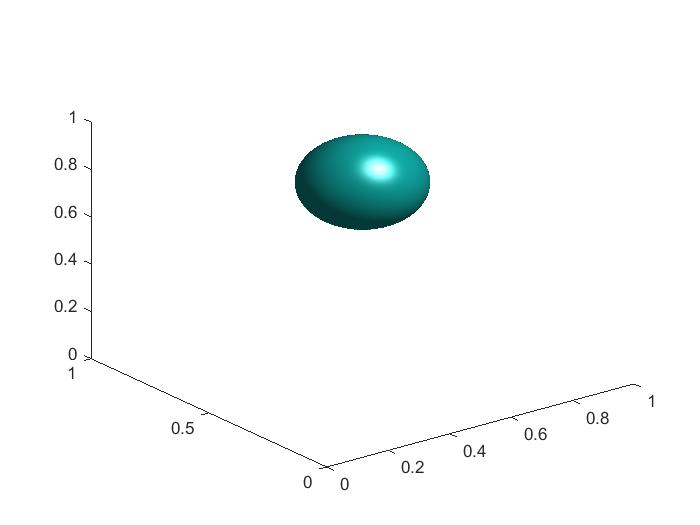}  \\ 
(a) Slice image of the true $n$ & (b) 3D image of the true $n$ \\ 
\includegraphics[width=3.5cm]{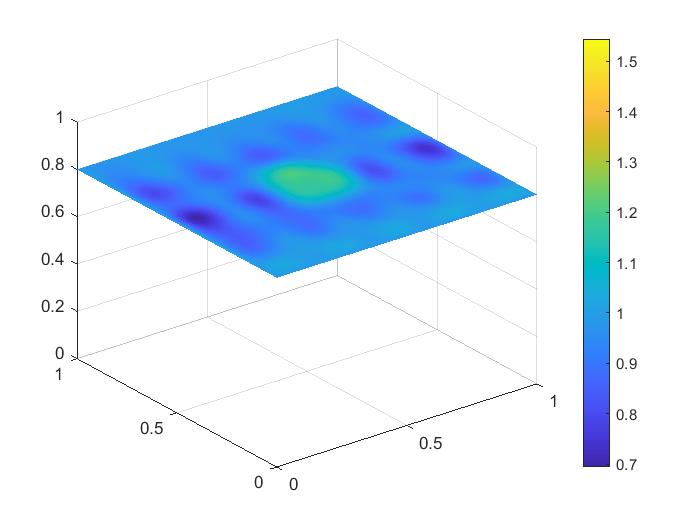}
& %
\includegraphics[width=3.5cm]{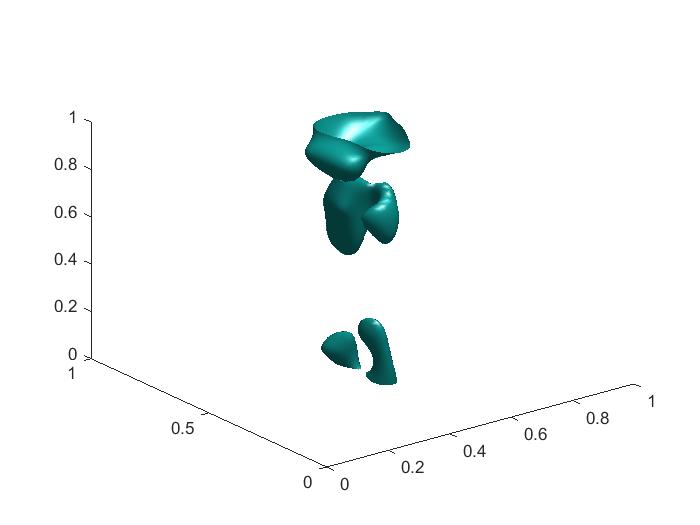}
\\ 
(c) Slice image of $n_{comp}$ for $\lambda=0$ & (d) 3D image of $n_{comp}$
for $\lambda=0$ \\ 
\includegraphics[width=3.5cm]{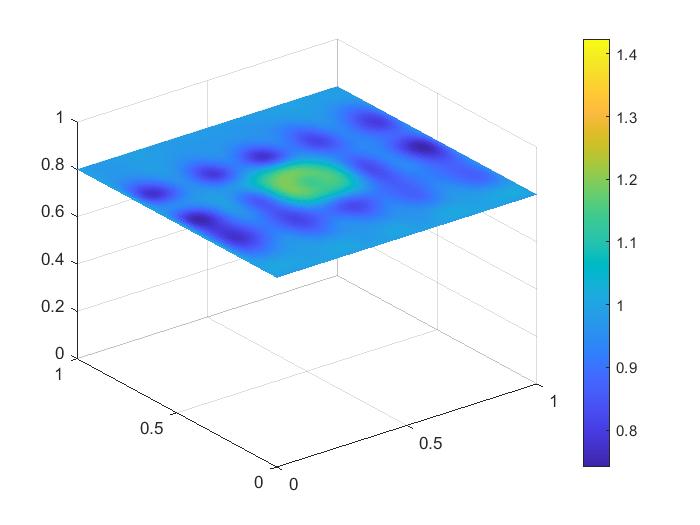}
& %
\includegraphics[width=3.5cm]{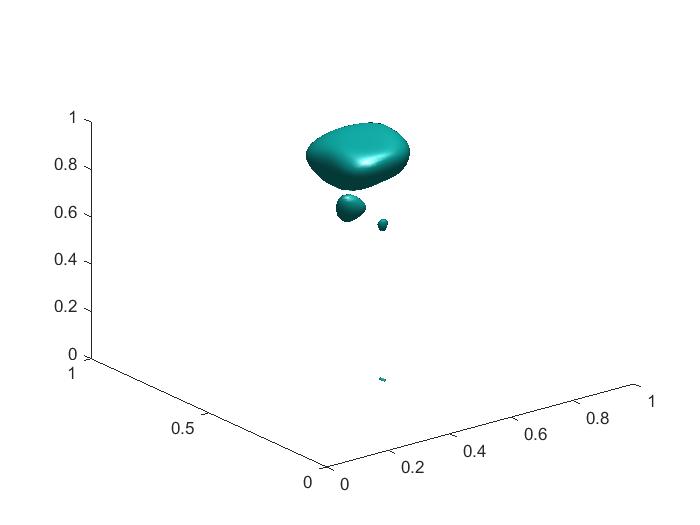}
\\ 
(e) Slice image of $n_{comp}$ for $\lambda=1$ & (f) 3D image of $n_{comp}$
for $\lambda=1$ \\ 
\includegraphics[width=3.5cm]{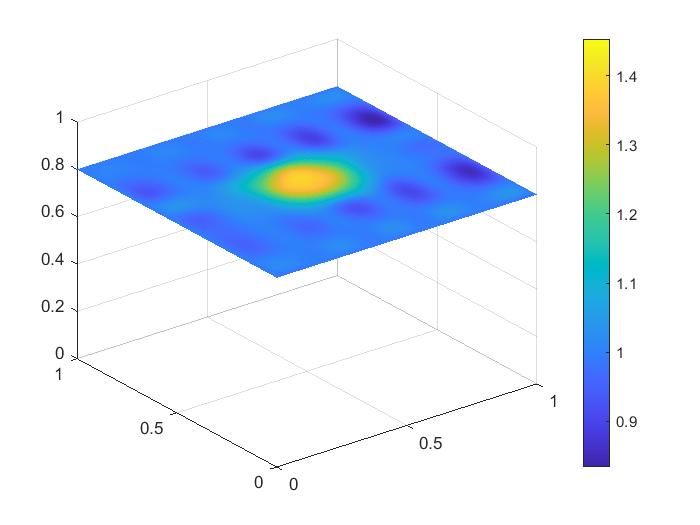}
& %
\includegraphics[width=3.5cm]{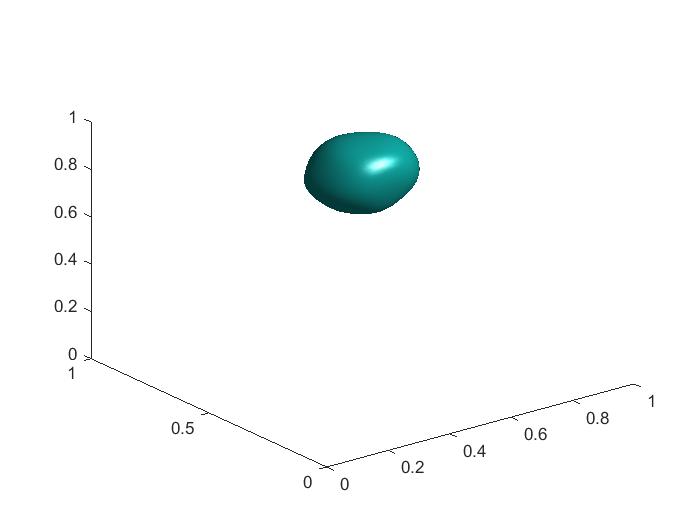}
\\ 
(g) Slice image of $n_{comp}$ for $\lambda=2$ & (h) 3D image of $n_{comp}$
for $\lambda=2$ \\ 
\includegraphics[width=3.5cm]{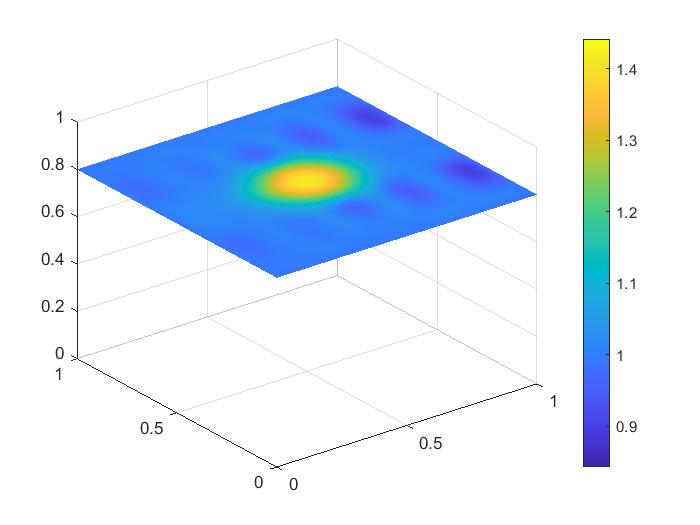}
& %
\includegraphics[width=3.5cm]{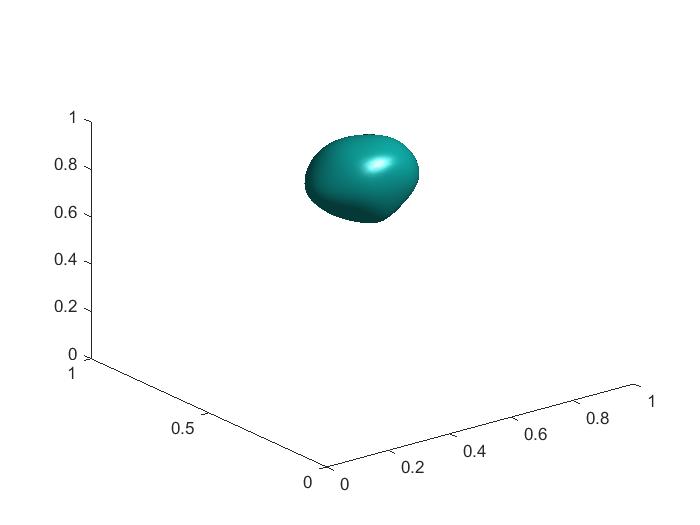}
 \\ 
(i) Slice image of $n_{comp}$ for $\lambda=3$ & (j) 3D image of $n_{comp}$
for $\lambda=3$  \\ 
\includegraphics[width=3.5cm]{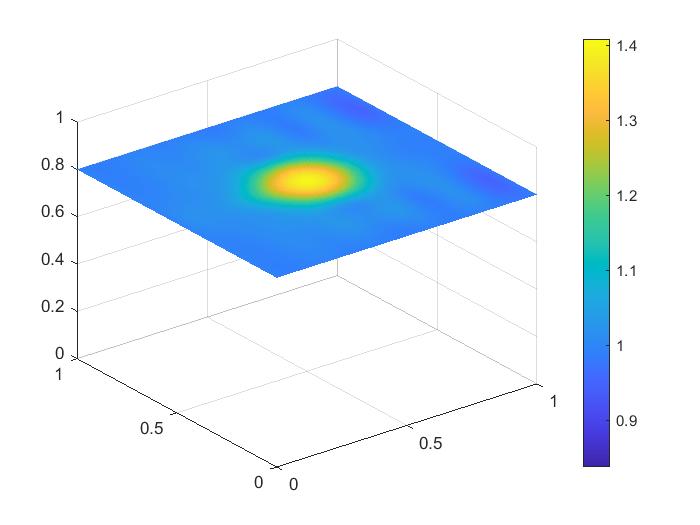}
& %
\includegraphics[width=3.5cm]{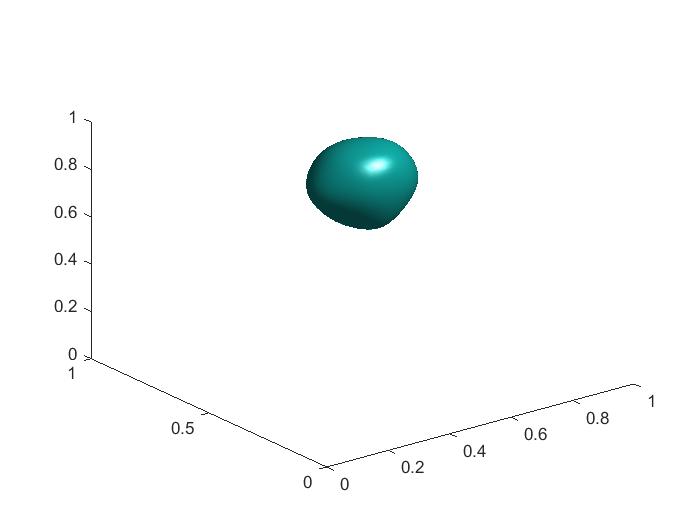}
\\ 
(k) Slice image of $n_{comp}$ for $\lambda=4$ & (l) 3D image of $n_{comp}$
for $\lambda=4$ 
\end{tabular}%
\end{center}
\caption{\emph{Results of Test 1. We test the effects of the parameter $%
\protect\lambda$. In this test, we fix $N=6$ and $\protect\lambda$ varies
from 0 to 4.}}
\label{example1}
\end{figure}

\begin{figure}[tbp]
\begin{center}
\begin{tabular}{cc}
\includegraphics[width=4cm]{Figures/3D-true-1in558-1-5-slice.jpg} & %
\includegraphics[width=4cm]{Figures/3D-true-1in558-1-5-sur.jpg} \\ 
(a)Slice image of the true $n$ & (b) 3D image of the true $n$ \\ 
\includegraphics[width=4cm]{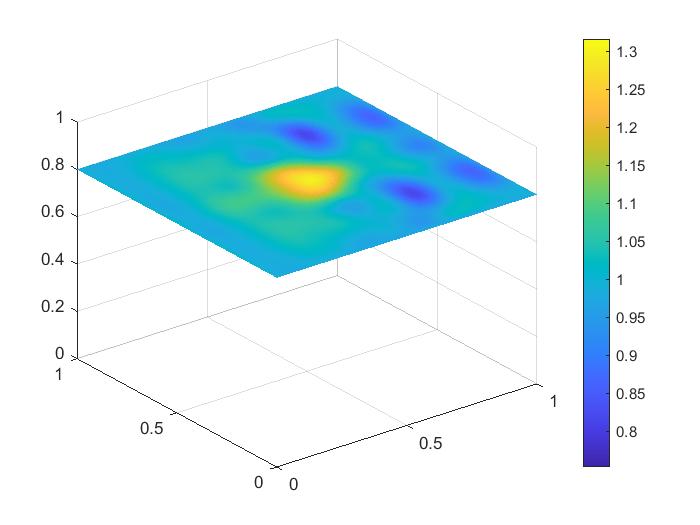}
& %
\includegraphics[width=4cm]{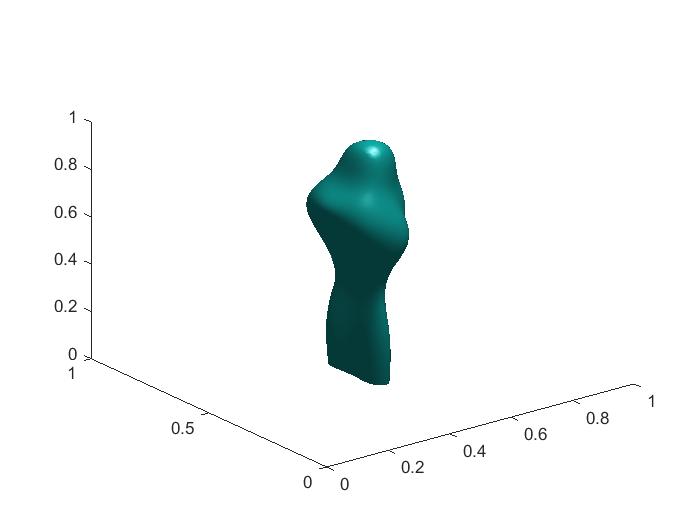}
\\ 
(c) Slice image of $n_{comp}$ for $N=4$ & (d) 3D image of $n_{comp}$ for $%
N=4 $ \\ 
\includegraphics[width=4cm]{Figures/1in-1-5-thick-015-A-2-A1-1-n-20-lam-4-beta-1e-4-N6-slice.jpg}
& %
\includegraphics[width=4cm]{Figures/1in-1-5-thick-015-A-2-A1-1-n-20-lam-4-beta-1e-4-N6-sur.jpg}
\\ 
(e) Slice image of $n_{comp}$ for $N=6$ & (f) 3D image of $n_{comp}$ for $%
N=6 $ \\ 
\includegraphics[width=4cm]{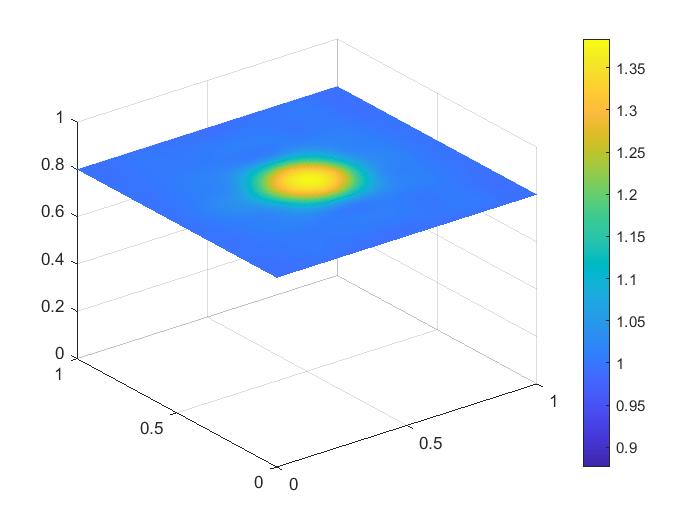}
& %
\includegraphics[width=4cm]{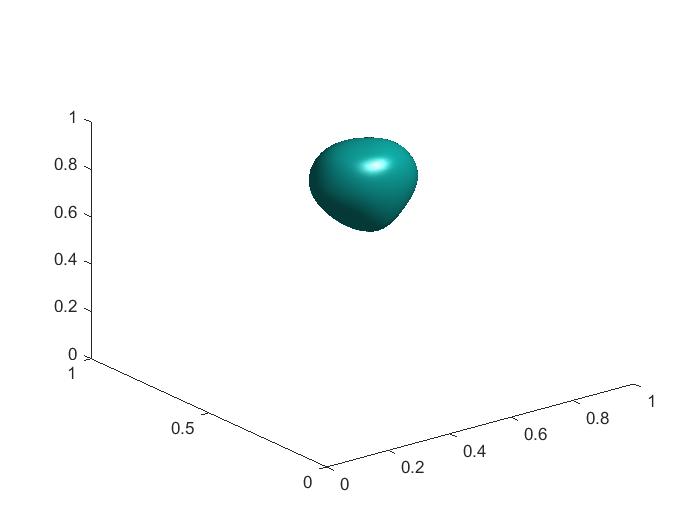}
\\ 
(g) Slice image of $n_{comp}$ for $N=8$ & (h) 3D image of $n_{comp}$ for $%
N=8 $ \\ 
\includegraphics[width=4cm]{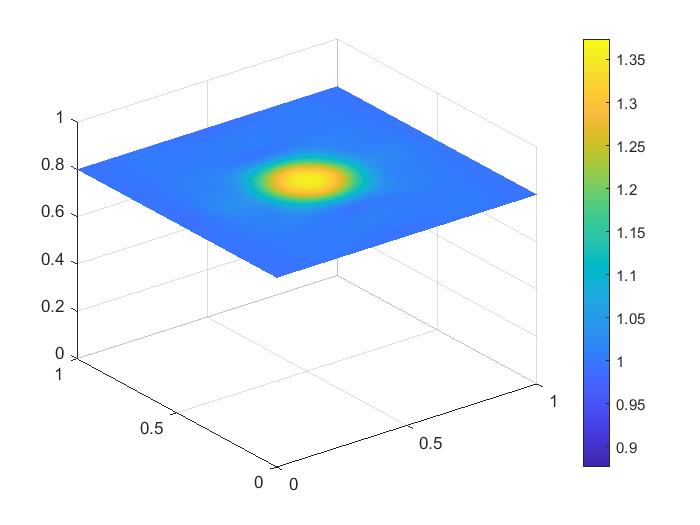}
& %
\includegraphics[width=4cm]{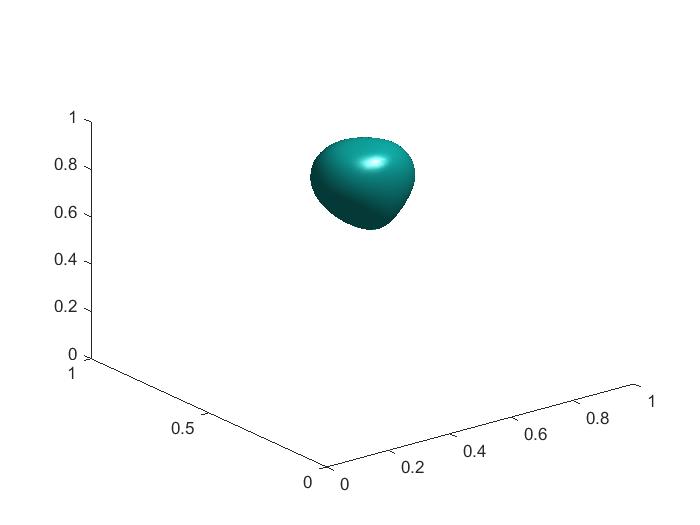}
\\ 
(i) Slice image of $n_{comp}$ for $N=10$ & (j) 3D image of $n_{comp}$ for $%
N=10$%
\end{tabular}%
\end{center}
\caption{\emph{Results of Test 2. We test the effects of the parameter $N$.
Here we fix the parameter $\protect\lambda=4$ and $N$ varies from 4 to 10.}}
\label{example2}
\end{figure}

\begin{figure}[tbp]
\begin{center}
\begin{tabular}{cc}
\includegraphics[width=4cm]{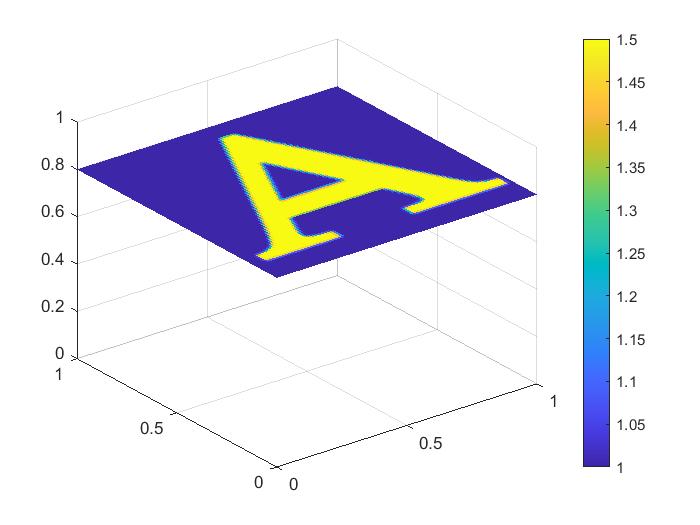} & %
\includegraphics[width=4cm]{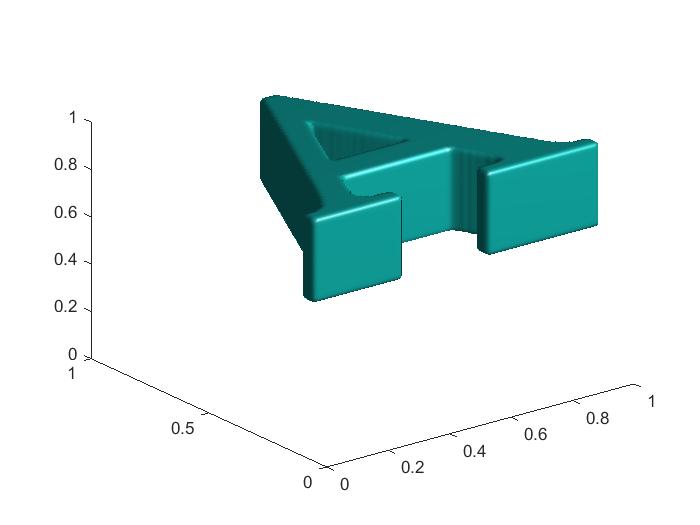}\\ 
(a) Slice image of the true $n$ & (b) 3D image of the true $n$ \\ 
\includegraphics[width=4cm]{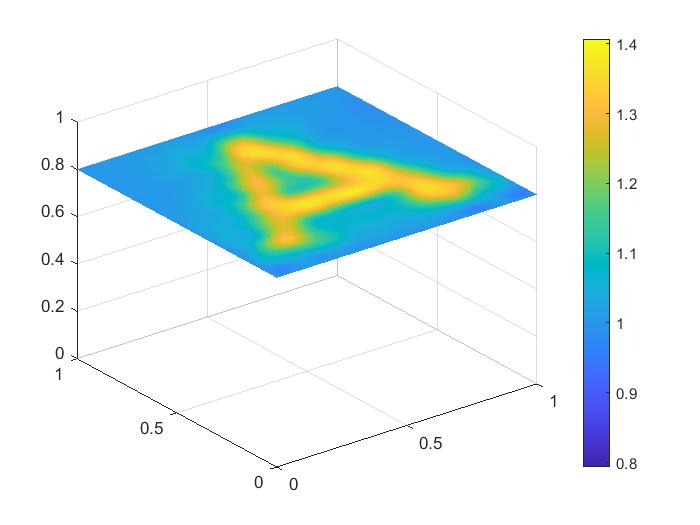}
& %
\includegraphics[width=4cm]{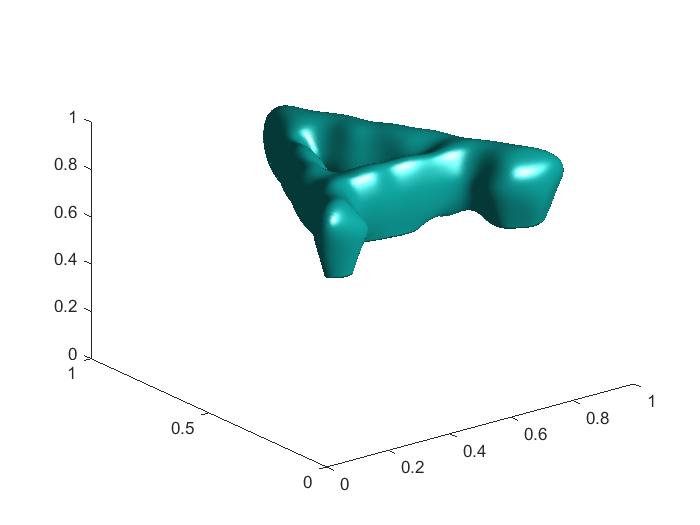}
 \\ 
(c) Slice image of $n_{comp}$ & (d) 3D image of $n_{comp}$
\end{tabular}%
\end{center}
\caption{\emph{Results of Test 3. Imaging of `A' shaped $n$ with $n=1.5$ in
it and $n=1$ outside. In this example, we set the parameter $N=6$ and $%
\protect\lambda=4$, see \eqref{5.10}. }}
\label{example3}
\end{figure}

\begin{figure}[tbp]
\begin{center}
\begin{tabular}{cc}
\includegraphics[width=4cm]{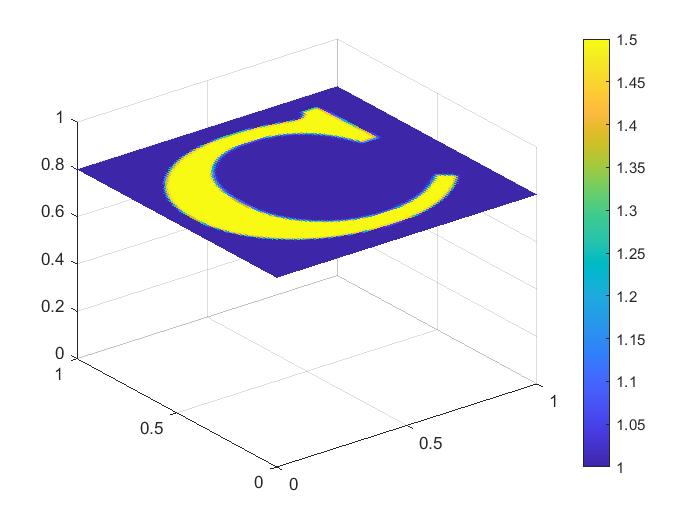} & %
\includegraphics[width=4cm]{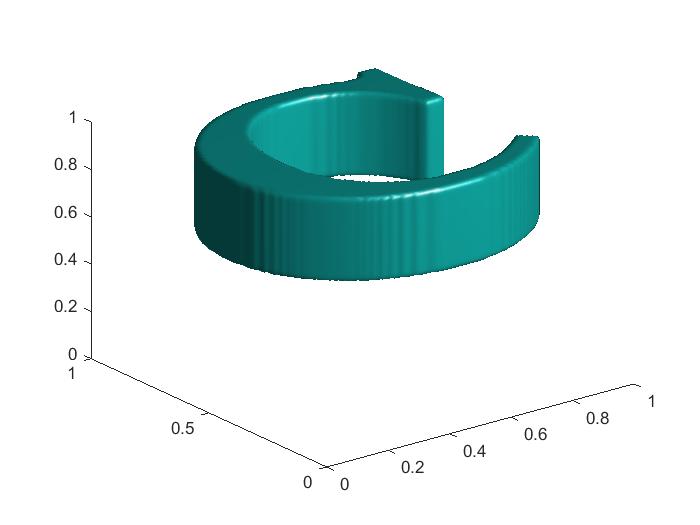} \\ 
(a) Slice image of the true $n$ & (b) 3D image of the true $n$ \\ 
\includegraphics[width=4cm]{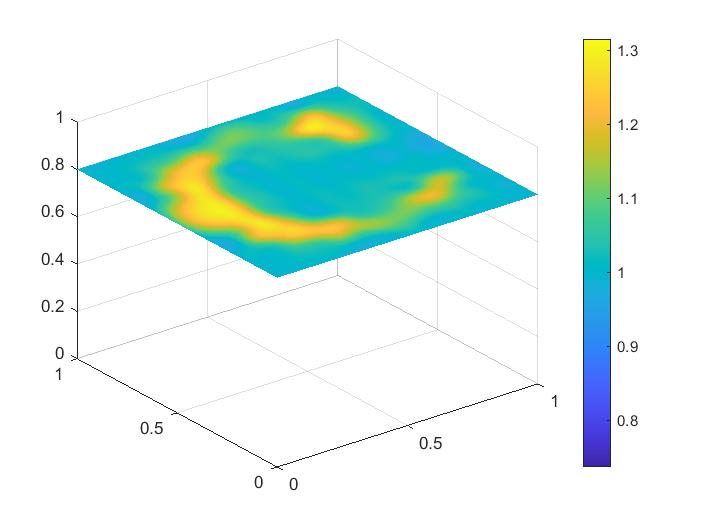}
& %
\includegraphics[width=4cm]{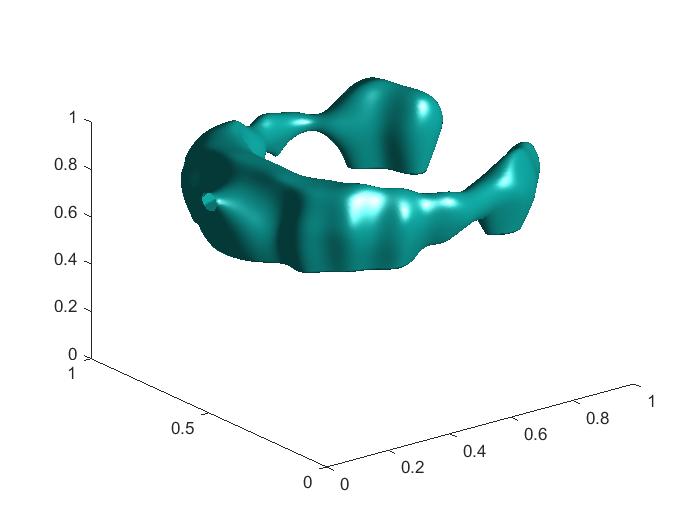}
\\ 
(c) Slice image of $n_{comp}$ & (d) 3D image of $n_{comp}$ 
\end{tabular}%
\end{center}
\caption{\emph{Results of Test 4. Imaging of `C' shaped $n$ with $n=1.5$ in
it and $n=1$ outside. In this example, we set the parameter $N=6$ and $%
\protect\lambda=4$, see \eqref{5.10}. }}
\label{example4}
\end{figure}

\begin{figure}[tbp]
\begin{center}
\begin{tabular}{cc}
\includegraphics[width=4cm]{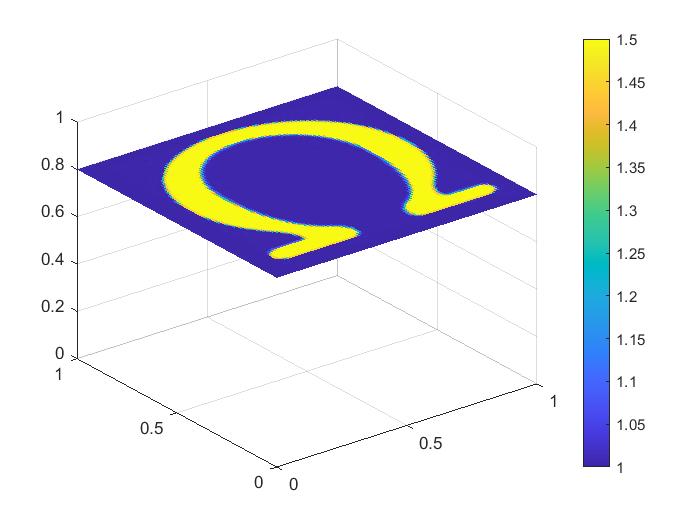}
& \includegraphics[width=4cm]{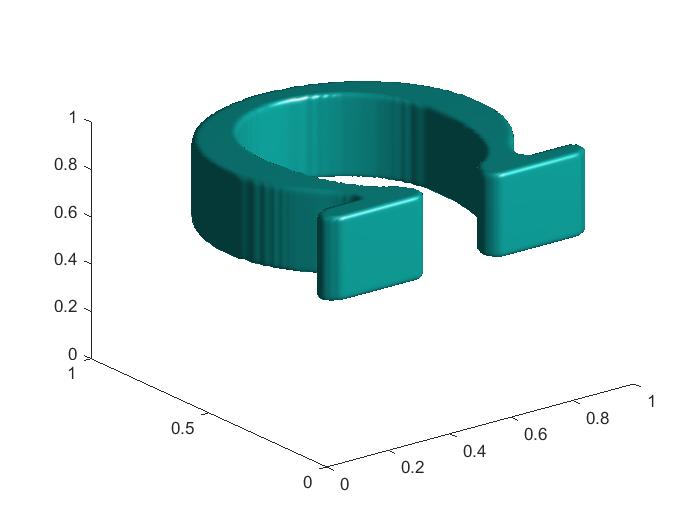}
\\ 
(a) Slice image of the true $n$ & (b) 3D image of the true $n$ \\ 
\includegraphics[width=4cm]{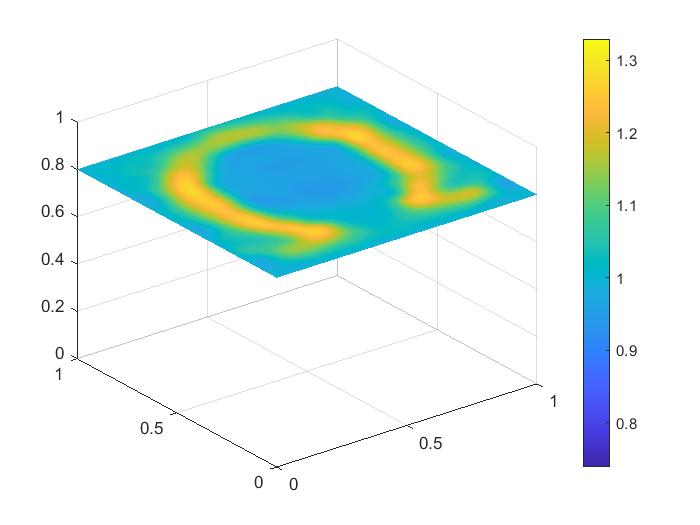}
& %
\includegraphics[width=4cm]{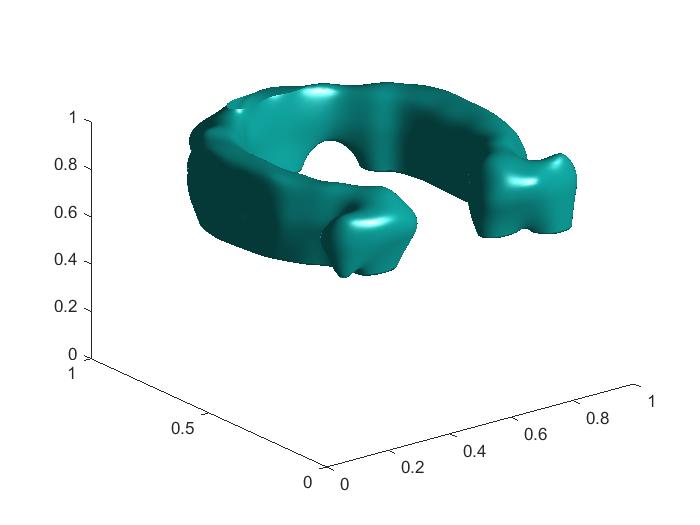}
\\ 
(c) Slice image of $n_{comp}$ & (d) 3D image of $n_{comp}$ 
\end{tabular}%
\end{center}
\caption{\emph{Results of Test 5. Imaging of ` $\Omega$' shaped $n$ with $%
n=1.5$ in it and $n=1$ outside. In this example, we set the parameter $N=6$
and $\protect\lambda=4$, see \eqref{5.10}. }}
\label{example5}
\end{figure}

The true inclusion/background contrast in the target function $m\left( 
\mathbf{x}\right) $ is 2.25:1 in all tests. On the other hand, it follows
from Figures \ref{example1}-\ref{example5} that the computed contrasts vary
between $1.3^{2}:1=1.69:1$ and $1.4^{2}:1=1.96:1.$ Also, some other
refinements of our results are desirable. We hope to obtain them in the
future.

\section*{CRediT Authorship Contribution Statement}

\textbf{1. Michael V. Klibanov:} Conceptualization, Methodology, Formal
Analysis, Supervision, Writing -- original draft, Writing -- review \&
editing.

\textbf{2. Jingzhi Li}: Conceptualization, Methodology, Supervision, Writing
-- original draft, Writing -- review \& editing.

\textbf{3. Wenlong Zhang:} Conceptualization, Investigation, Software,
Validation, Visualization, Writing -- review \& editing

\section{\noindent Declaration of competing interest}

\noindent The authors declare that they have no known competing financial
interests or personal relationships that could have appeared to influence
the work reported in this paper.

\section*{Acknowledgement}

\noindent The research of W. Zhang is supported by the National Natural
Science Foundation of China No. 11901282 and the Shenzhen Sci-Tech Fund No.
RCBS20200714114941241.


\end{document}